\documentclass[12pt]{article}

\usepackage{amsfonts}
\usepackage{epsfig}
\usepackage[percent]{overpic}
\usepackage{amsmath}
\numberwithin{equation}{section}

\usepackage[letterpaper,margin=0.85in]{geometry}

\usepackage[utf8]{inputenc}
\usepackage{amssymb}
\usepackage{graphics}
\usepackage{verbatim}
\usepackage{bbm}
\usepackage{color}
\usepackage{amsthm}
\usepackage[all]{xy}
\usepackage{makeidx}
\usepackage{enumerate}
\usepackage{mathrsfs}

\let\oldbibliography\thebibliography
\renewcommand{\thebibliography}[1]{%
  \oldbibliography{#1}%
  \setlength{\itemsep}{-1.2mm}%
}

\theoremstyle{plain}
\newtheorem{thm}{Theorem}[section]

\newtheorem{lem}[thm]{Lemma}
\newtheorem{prop}[thm]{Proposition}

\theoremstyle{definition}
\newtheorem{defn}[thm]{Definition}

\newtheoremstyle{myremark}% name
  {3pt}%      Space above
  {3pt}%      Space below
  {\small \rmfamily}%         Body font
  {5pt}%         Indent amount (empty = no indent, \parindent = para indent)
  {\rmfamily}% Thm head font
  {:}%        Punctuation after thm head
  {.5em}%     Space after thm head: " " = normal interword space;
  {}%         Thm head spec (can be left empty, meaning `normal')

\theoremstyle{myremark}
\newtheorem*{remark}{\textit{Remark}}

\def\R{\mathbb{R}}

\def\N{\mathbb{N}}

\def\cL{\mathcal{L}}

\def\cU{\mathcal{U}}
\def\cV{\mathcal{V}}

\def\cX{\mathcal{X}}

\def\Pp{\mathsf{P}}
\def\eps{\epsilon} 

\def\txtc{{\textnormal{c}}}
\def\txtd{{\textnormal{d}}}

\def\txti{{\textnormal{i}}}

\def\txtD{{\textnormal{D}}}

\def\Leb{{\textnormal{Leb}}}

\def\ra{\rightarrow}
\def\I{\infty}

\newcommand{\bpf}[1][Proof]{{\noindent {\sc #1: }}}
\newcommand{\epf}{{{\hfill $\Box$ \smallskip}}}

\newcommand{\be}{\begin{equation}}
\newcommand{\ee}{\end{equation}}
\newcommand{\benn}{\begin{equation*}}
\newcommand{\eenn}{\end{equation*}}
\newcommand{\bea}{\begin{eqnarray}}
\newcommand{\eea}{\end{eqnarray}}
\newcommand{\beann}{\begin{eqnarray*}}
\newcommand{\eeann}{\end{eqnarray*}}

\newcommand{\myendex}{$\blacklozenge$\end{ex}}
\newcommand{\myendexerc}{$\lozenge$\end{exerc}}
\newcommand{\myendpexerc}{$\lozenge$\end{pexerc}}

\def\XXint#1#2#3{{\setbox0=\hbox{$#1{#2#3}{\int}$}
\vcenter{\hbox{$#2#3$}}\kern-.5\wd0}}

\begin{document}

\author{Tobias Hurth\thanks{Institut de Math{\'e}matiques, Universit{\'e} de Neuch{\^a}tel, 2000 Neuch{\^a}tel, Switzerland}~~and~Christian Kuehn\thanks{Technical 
University of Munich, Faculty of Mathematics, 85748 Garching bei M\"unchen, 
Germany}}
 
\title{Random Switching near Bifurcations}

\maketitle

\begin{abstract}
The interplay between bifurcations and random switching processes of vector fields is studied. More precisely, we provide a classification of piecewise deterministic Markov processes arising from stochastic switching dynamics near fold, Hopf, transcritical and pitchfork bifurcations. We prove the existence of invariant measures for different switching rates. We also study, when the invariant measures are unique, when multiple measures occur, when measures have smooth densities, and under which conditions finite-time blow-up occurs. We demonstrate the applicability of our results for three nonlinear models arising in applications. 
\end{abstract}

%%%%%%%%%%%%%%%%%%%%%%%%%%%%%%%%%%%%%%%%%%%%%%%%%%%%%%%%%%%%%%%%%%%%%%%%%%%%%%%%
\section{Introduction}
\label{sec:intro}

In this work we study the dynamics of randomly switched ordinary differential
equations (ODEs) of the form
\be
\label{eq:ODEintro}
\frac{\txtd x}{\txtd t}=x'=f(x,p),\qquad x=x(t)\in\R^d,~x(0)=:x_0,
\ee
near bifurcation points. More precisely, we select two parameters
$p=p_\pm\in\R$ so that~\eqref{eq:ODEintro} has non-equivalent dynamics~\cite{Kuznetsov}, 
which are separated by a distinguished bifurcation point $p_*\in(p_-,p_+)$. Then
we look at the piecewise-deterministic Markov process (PDMP) generated by
switching between the vector fields $f(x,p_-)$ and $f(x,p_+)$. This idea
is motivated by several observations. Here we just name a few:

\begin{enumerate}
 \item[(M1)] In parametric families of vector fields, bifurcations occur generically. Therefore, they are immediately relevant for the study of PDMPs as well. In addition, the interplay between random switching and bifurcation points is not studied well enough yet.
 \item[(M2)] From the perspective of PDMPs, this setting provides natural examples to test and extend the general theory of invariant measures.
 \item[(M3)] Stochastic bifurcation theory is a very active area, where still many questions remain open. Hence, studying a well-defined set of standard cases involving bifurcations and stochasticity is highly desirable.
 \item[(M4)] Parameters in many models are usually only known via a possible distribution and not exactly. Therefore, our work contributes to the uncertainty quantification for nonlinear systems arising in applications. 
\end{enumerate}

Before describing our main results, we briefly review some of the background from PDMPs and from nonlinear dynamics to provide a broader perspective.\medskip 

The study of randomly switched deterministic vector fields goes back at least to the works of Goldstein \cite{Goldstein} and Kac \cite{Kac}. The set-up can informally be described as follows: Given a starting point $x_0 \in \R^d$ and an initial vector field $f_i$ taken from a finite collection $\{f_j\}$ of smooth vector fields on $\R^d$, we follow the flow along $f_i$ starting at $x_0$ for an exponentially distributed random time. Then a switch occurs, meaning that $f_i$ is replaced with a new vector field $f_j$, $j\neq i$. We flow along $f_j$ for another exponential time and switch again. This yields a continuous and piecewise smooth trajectory in $\R^d$ that is, however, not the trajectory of a Markov process. To obtain a Markov process, one needs to supplement the switching process on $\R^d$ with a second stochastic process that keeps track of the driving vector field. The resulting two-component process belongs to the class of piecewise deterministic Markov processes (PDMPs).  

PDMPs were first introduced by Davis~\cite{Davis_article} in an even more general setting. For instance, PDMPs may involve jumps not only on the collection of vector fields but also on $\R^d$ ~\cite{Davis,Malrieu_2015}. The class of PDMPs considered in this article is also known under the names of hybrid systems \cite{Yin} and random evolutions \cite{Hersh},~\cite[Chapter~12]{Ethier_Kurtz}. Randomly switched vector fields have applications to areas such as ecology \cite{Lobry}, gene regulation \cite{Bressloff}, molecular motors \cite{Gabrielli}, epidemiology \cite{Cui}, queueing theory \cite{Anisimov}, and climate science \cite{Majda}, to name just a few. 

Aside from their uses in modeling, randomly switched vector fields have intriguing theoretical properties. For example, switching between stable vector fields can result in an unstable situation, and vice versa. Recently, examples of randomly switched vector fields were found that exhibit such a reversal of stability for almost all realizations of switching times~\cite{Zitt,Lawley}. Another interesting phenomenon is the regularizing effect random switching can have on a dynamical system. For example, random switching between two Lorenz vector fields with just slightly different parameter values induces an invariant probability measure that is absolutely continuous with respect to Lebesgue measure on $\R^3$, whereas the dynamics associated to each individual vector field concentrate on attractors of Lebesgue measure zero~\cite{Strickler,Bakhtin}. Another recent topic is the ergodic theory for randomly switched vector fields. Important contributions to the question whether a switching system on a noncompact state space admits an invariant probability measure were made in~\cite{Le_Borgne,B17,BS17}. In~\cite{Bakhtin,Benaim}, it was shown that a H\"ormander-type hypoellipticity condition on the vector fields at an accessible point yields uniqueness and absolute continuity of the invariant probability measure. 

In this work we focus on the invariant probability measure aspect and relate it to bifurcation points. Bifurcation theory~\cite{GH,Kuznetsov} has become one of the most widely used techniques to study nonlinear systems~\cite{Strogatz}. Informally, the main idea is to study vector fields under parameter variation and to determine at which points the dynamics changes fundamentally, i.e., to detect the points where the phase portraits of the vector fields are not topologically equivalent upon small parameter variation. Almost full classification results exist for bifurcations with relatively few parameters, i.e., codimension one or two. These results provide suitable unfoldings, which are basically partitions of parameter space into non-equivalent phase portraits~\cite{Kuznetsov}.

Recently, substantial interest has been focused on understanding the interplay between stochasticity and bifurcations. Yet, the setting in almost all of these works is focused on either stochastic differential equations (SDEs) involving (space-)time stochastic forcing processes~\cite{ArnoldSDE,BerglundGentz}, or less frequently on random differential equations (RDEs) with a fixed random parameter distribution~\cite{BredenKuehn,NairSarkarSujith}. Particularly interesting dynamics seems to appear for SDEs in oscillatory situations~\cite{BaudelBerglund,EngelLambRasmussen,SadhuKuehn}. Recently numerical and semi-analytical work shows that interesting effects also occur for switched systems near bifurcations~\cite{KuehnQuenched}. Therefore, it is very natural that one should try to link PDMPs with bifurcation theory.\medskip 

In this paper, we provide a full mathematical classification of the PDMPs associated to~\eqref{eq:ODEintro} switched near local bifurcations for codimension one bifurcations. We not only include the generic fold and Hopf bifurcations but also study the frequently occurring one-parameter transcritical and pitchfork bifurcations. We prove under which conditions on the switching rates invariant measures occur, when they are unique, when their densities are smooth, and we also provide explicit formulas for these densities in certain cases. In addition, we prove finite-time blow-up results for certain parameter regimes. In summary, our theorems provide building blocks, which can be employed in various PDMPs. In addition, we demonstrate that we may also derive insights from our results in three nonlinear models arising respectively in ecology, nonlinear oscillations, and collective motion.\medskip

The paper is structured as follows: In Section~\ref{sec:background} we provide more technical background from local bifurcation theory and PDMPs. In Section~\ref{sec:stable} we focus on all cases where below and above the bifurcation point there are non-trivial trapping regions. In these cases we characterize the occurring invariant probability measures completely. In Section~\ref{sec:unstable} we consider the cases with only one non-trivial trapping region. We again study the invariant measures in full detail but now also finite-time blow-up can appear. In Section~\ref{sec:applications}, we indicate how our results can be used in three models arising from applications.

%%%%%%%%%%%%%%%%%%%%%%%%%%%%%%%%%%%%%%%%%%%%%%%%%%%%%%%%%%%%%%%%%%%%%%%%%%%%%%%%
\section{Background}
\label{sec:background}

We briefly recall the technical background needed from the two main areas we consider
in this work. Hence, this section mainly serves as a reference and to fix the
notation. Readers familiar with local bifurcation theory~\cite{GH,Kuznetsov} 
and PDMPs can skip ahead to Section~\ref{sec:stable}.

%%%%%%%%%%%%%%%%%%%%%%%%%%%%%%%%%%%%%%%%%%%%%%%%%%%%%
\subsection{Local Bifurcation Theory}
\label{sec:biftheory}

Consider an ordinary differential equation (ODE) given by
\be
\label{eq:ODE}
\frac{\txtd x}{\txtd t}=x'=f(x,p),\qquad x=x(t)\in\R^d,~x(0)=:x_0,
\ee
where $p\in\R$ is the (main) bifurcation parameter, and we assume that 
the vector field
$f:\R^d\times \R\ra \R^d$ is sufficiently smooth; in particular, in what
follows $f\in C^3(\R^d \times \R,\R^d)$ is going to suffice. We also refer to 
$\R^d$ as the phase space of~\eqref{eq:ODE}. The phase space together with 
the foliation by trajectories $x(t)$ is called phase portrait. Suppose $x_*$ is an 
equilibrium point (or steady state) of~\eqref{eq:ODE} for the parameter 
value $p_*$ so that $f(x_*,p_*)=0$. Without loss of generality, upon
translating coordinates in the phase space $\R^d$ and the parameter space
$\R$, we may assume that $(x_*,p_*)=(0,0,\ldots,0)=:0$. Consider the 
linearized problem near the steady state
\be
\label{eq:linODE}
X'=\txtD_xf(0)X=AX,\qquad X=X(t)\in\R^d.
\ee
Then $x_*$ is called hyperbolic if the matrix $A\in\R^{d\times d}$
has no spectrum on the imaginary axis. In the hyperbolic case, the 
Hartman-Grobman Theorem~(see e.g.~\cite{Teschl}) implies that the systems~\eqref{eq:ODE} 
and \eqref{eq:linODE} are locally topologically equivalent, i.e., small 
parameter variations for $p\in(-p_0,p_0)$, $p_0>0$, do not qualitatively alter 
the phase portrait as the hyperbolic structure of $A$ is robust under 
small parameter perturbations. More precisely, for any two parameter
values $p_1,p_2\in(-p_0,p_0)$, there exists a homeomorphism $h:\R^d\ra \R^d$
such that the phase portraits of $f(x,p_1)$ and $f(x,p_2)$ are mapped 
to each other by $h$ preserving the direction of time on trajectories. 

Suppose $A$ is not hyperbolic so that 
$\textnormal{spec}(A)\cap \txti \R\neq \emptyset$. A local bifurcation 
occurs at $p_*=0$ if for any $p_0>0$ and any open neighbourhood $\cU=\cU(0)$ 
of $x_*=0$, there exist two locally (wrt $\cU$) non-homeomorphic phase 
portraits of~\eqref{eq:ODE} for two values $p_1,p_2\in(-p_0,p_0)$. In 
particular, a bifurcation just corresponds to the appearance of a topologically 
non-equivalent phase portrait under parameter variation. 

The general strategy to analyze bifurcation problems~\cite{Kuznetsov,GH} 
proceeds as follows:
(I) the system is reduced to the dimension $d_\txtc$ of $\ker(A)$ using
a center manifold $W^\txtc_{\textnormal{loc}}(0)$, (II) on 
$W^\txtc_{\textnormal{loc}}(0)$ one uses smoothness to Taylor-expand
the vector field and then simplifies it using coordinate changes, and
(III) one proves that a finite number of polynomial terms is locally
sufficient to determine the topological equivalence class so a truncation
yields a finite-degree polynomial vector field.
The steps (I)-(III) lead to different classes/families of polynomial 
vector fields, also called normal forms, depending upon degeneracy 
of $\textnormal{spec}(A)$ and depending upon a finite number of partial 
derivatives of $f$. 

\begin{figure}[htbp]
	\centering
	\begin{overpic}[width=0.4\textwidth]{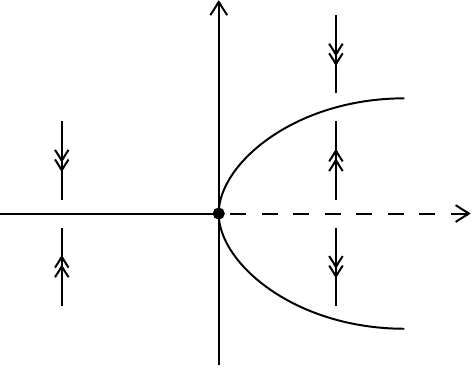}
		\put(95,37){\scalebox{1}{$p$}}
		\put(50,73){\scalebox{1}{$x$}}
		\put(22,15){\scalebox{1}{$S_{p \leq 0}$}}
		\put(80,15){\scalebox{1}{$S_{p > 0}$}}
		\end{overpic}
	\caption{\label{fig:01}Sketch of the bifurcation diagram for 
	the supercritical pitchfork bifurcation normal form~\eqref{eq:pitchforkex}.
	There are two classes of topologically non-equivalent phase portraits
	here denoted by $S_{p \leq 0}$ and $S_{p >0}$.}
\end{figure}
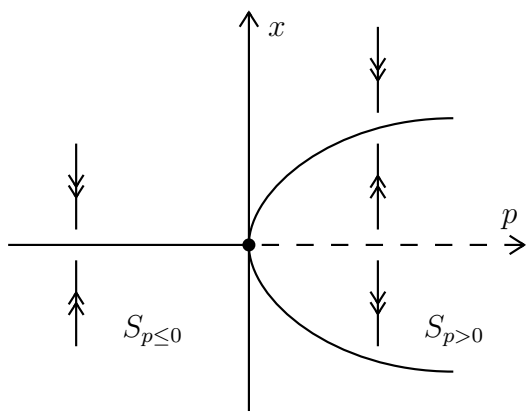

In this work we shall focus on the four most common bifurcations
used in practical applications for $d_\txtc=1$ and $d_\txtc=2$, which
just require a single bifurcation parameter $p$, and where the system 
has already been reduced to normal form. These cases will be the 
fold, Hopf, transcritical, and pitchfork bifurcations. As a motivating 
example, consider the supercritical pitchfork normal form
\be
\label{eq:pitchforkex}
x'=px-x^3,\qquad x\in\R,~p\in\R.
\ee
Clearly, the equilibrium $x_*=0$ undergoes a bifurcation at $p_*=0$
as the phase portrait for $p<0$ has one globally stable equilibrium,
while the phase portrait for $p>0$ has three equilibria. For $p>0$,
we find that $x_*=0$ is unstable while the equilibria 
$x_\pm=\pm\sqrt{p}$ are both locally stable; see also Figure~\ref{fig:01}.

However, note that from the viewpoint of applications, treating $p$ just
as a static parameter is not always realistic. This approach presumes $p$ is
changed infinitely slowly to bring the system to and across the bifurcation
point. One option is to consider the case when $p$ is just switched across 
the bifurcation point (e.g., consider shot noise
effects, control action, activation of external interfaces of the system, etc).
As an example, consider the problem of switching between $p<0$ and $p>0$ 
in the context of the pitchfork normal form~\eqref{eq:pitchforkex}. This leads us 
naturally to consider piecewise deterministic Markov processes as introduced in the 
next section.

%%%%%%%%%%%%%%%%%%%%%%%%%%%%%%%%%%%%%%%%%%%%%%%%%%%%%
\subsection{Piecewise Deterministic Markov Processes}
\label{ssec:PDMP}

In this subsection we introduce a class of PDMPs characterized by Poissonian random switching between a finite number of deterministic vector fields.  Let $I$ be a finite index set, and let $(f_i)_{i \in I}$ be a collection of vector fields on $\R^d$ with some degree of smoothness.  To introduce the basic framework, we just assume that $(f_i)_{i \in I}$ are in $C^1(\R^d, \R^d)$, but for some of the results stated below higher degrees of smoothness are required.  To be able to associate flows to the vector fields, we assume in addition that $(f_i)_{i \in I}$ are forward complete, i.e. for any $x_0 \in \R^d$ the initial-value problem 
$$ 
x' = f_i(x), \quad x(0) = x_0
$$ 
has a unique solution $t \mapsto \Phi_i^t(x_0)$ that is defined for all $t \geq 0$.  Given a starting point $x_0 \in \R^d$ and an initial vector field $f_i$, the random dynamical system we consider follows the flow associated to $x_0$ and $f_i$ for a random time. Then a switch occurs, which means that the driving vector field $f_i$ is replaced by a new vector field $f_j$ chosen at random from $\{f_k: k \in I \setminus \{i\}\}$. Again, the system flows along $f_j$ for a random time until another switch occurs, etc.  The stochastic process $X = (X_t)_{t \geq 0}$ that records the position of the switching trajectory on $\R^d$ is not Markov because knowing $(X_s)_{s \leq t}$ lets us infer the driving vector field at time $t$.  If the times between consecutive switches are exponentially distributed and independent conditioned on the sequence of driving vector fields, and if the vector fields are chosen according to a Markov chain on $I$, then the two-component process $(X,E)$ is already Markov, where $E_t \in I$ gives the index of the driving vector field at time $t$. For more general distributions of switching times, one needs to adjoin a third component that keeps track of the time elapsed since the latest switch.  It is possible to consider the situation where the rate of switching depends continuously on the location of the switching trajectory on $\R^d$ ~\cite{Benaim},~\cite{Faggionato}.  For simplicity we assume that the switching rates do not depend on the process $X$. We can then give the following rigorous description of $(X,E)$. Let $E = (E_t)_{t \geq 0}$ be an irreducible continuous-time Markov chain on the state space $I$.  Let $X = (X_t)_{t \geq 0}$ be the solution to the control problem 
$$ 
X_t = x + \int_0^t f_{E_s}(X_s) \ ds. 
$$ 
The Markov process $(X,E)$ has infinitesimal generator $L$ acting on functions $g: \R^d \times I \to \R$ that are smooth in $x$ according to 
\begin{equation}    \label{eq:generator}
Lg(x,i) = \langle f_i(x), \nabla_x g(x,i) \rangle + \sum_{j \neq i} \lambda _{i,j} (g(x,j) - g(x,i)),  
\end{equation}
where $\lambda_{i,j}$ is the rate at which $E$ transitions from state $i$ to state $j$. 
We denote the Markov semigroup of $(X,E)$ by $(\Pp^t)_{t \geq 0}$ or just by $(\Pp^t)$. An \emph{invariant probability measure} (IPM) of $(\Pp^t)$ is a probability measure $\mu$ on $\R^d \times I$ such that $\mu = \mu \Pp^t$ for all $t \geq 0$.  Below, we collect some results on existence, uniqueness and absolute continuity for IPM of $(\Pp^t)$ that have been established in the literature. 

We call a set $M \subset \R^d$ \emph{positive invariant} if $M$ is positive invariant under the flows $(\Phi_i)_{i \in I}$ associated with the vector fields $(f_i)_{i \in I}$, i.e. if for any $x \in M$, $i \in I$ and $t \geq 0$, we have $\Phi_i^t(x) \in M$. Thus, trajectories of $X$ starting in a positive invariant set $M$ or entering $M$ at some time remain in $M$ for all future times. If there is a compact positive invariant set $M$, existence of an IPM is guaranteed by the Krylov--Bogoliubov method~\cite[Theorem~3.1.1]{DaPrato}, which applies because $(X,E)$ is Feller~\cite[Proposition~2.1]{Benaim}. In the noncompact situation, an IPM is guaranteed to exist, for instance, if $(X,E)$ is on average contracting~\cite[Corollary~1.11]{Le_Borgne}.  By Harris's ergodic theorem, existence also holds if the semigroup $(\Pp^t)$ admits a Lyapunov function as well as a minorizing measure $\nu_K$ for every compact set $K\subset \R^d$.   

Recall that the \emph{support} of a Borel measure $\mu$ on $\R^d \times I$ is the set of points $(x,i) \in \R^d \times I$ such that $\mu(U \times \{i\}) > 0$ for every open neighborhood $U \subset \R^d$ of $x$.  If $x_* \in \R^d$ is an equilibrium for each of the vector fields $f_i$, then the product of the Dirac measure at $x_*$ and the unique IPM $\nu$ of the continuous-time Markov chain $E$ is a trivial IPM for $(\Pp^t)$.  If $M \subset \R^d$ is a compact positive invariant set containing such a common equilibrium $x_*$, the Krylov--Bogoliubov method is not sufficient to decide whether there are any additional IPM whose support is contained in $M \times I$.  This more subtle existence question can often be addressed using the theory of \emph{stochastic persistence} as developed by Bena\"im~\cite{B17}, and as applied to the case of a common equilibrium by Bena\"im and Strickler~\cite{Strickler}.  

We now outline an existence result from~\cite{Strickler} that will be needed later on. Let $M \subset \R^d$ be a compact positive invariant set containing the point $x_* = 0$, which we assume to be an equilibrium for all vector fields $f_i$. As a technical condition, we also require that there is $\delta > 0$ such that whenever $x \in M$ and $\|x\| \leq \delta$, then the entire line segment from $0$ to $x$ is contained in $M$. For $i \in I$, let 
$$ 
A_i = Df_i(0)
$$ 
be the Jacobian matrix of $f_i$ at $0$. Then, the cone 
$$ 
C_M = \{tx: t \geq 0, \ x \in M, \ \|x\| \leq \delta\} \subset \R^d
$$ 
is positive invariant with respect to the flows of the linear vector fields given by $(A_i)_{i \in I}$. On $C_M \times I$, we define the PDMP $(Y,E)$, which is obtained from $(X,E)$ by replacing each vector field $f_i$ with its linearization $A_i$. Whenever $Y_t \neq 0$, we define the angular process 
$$
\Theta_t = \frac{Y_t}{\|Y_t\|},  
$$ 
which evolves on the compact set $S^{d-1} \cap C_M$. By Krylov--Bogoliubov, $(\Theta,E)$ admits at least one IPM.  For any IPM $\nu$ of $(\Theta,E)$, define the \emph{average growth rate} as 
$$ 
\Lambda(\nu) = \sum_{i \in I} \int_{S^{d-1} \cap C_M} \theta^{\top} A_i \theta \ \nu(d \theta \times \{i\}). 
$$ 
Since $\|Y_t\|$ satisfies 
$$ 
\frac{d}{dt} \|Y_t \| = \Theta_t^{\top} A_{E_t} \Theta_t \| Y_t\|, 
$$ 
Birkhoff's ergodic theorem implies that for almost every realization of $(\Theta, E)$ with initial distribution $\nu$, we have 
$$ 
\lim_{t \to \infty} \frac{\ln(\|Y_t\|)}{t} = \Lambda(\nu). 
$$ 
Recall that an IPM $\nu$ of a Markov process with Markov semigroup $(\Pp^t)$ and state space $\cX$ is called \emph{ergodic} if $\nu(A) \in \{0,1\}$ for every measurable $A \subset \cX$ such that for all $t \geq 0$, $\Pp^t_x(A) = 1$ for $\nu$-almost every $x \in A$. Let $\Lambda^{-}$ denote the infimum and $\Lambda^+$ the supremum of $\Lambda(\nu)$ over all ergodic IPM $\nu$ of $(\Theta,E)$. In many situations of interest, $(\Theta,E)$ has exactly one IPM, so $\Lambda^- = \Lambda^+$. 

\begin{defn}
We call a point $x \in \R^d$ 
\begin{enumerate}
\item \emph{reachable from $y \in \R^d$} if there is a finite sequence of indices $i_1, \ldots, i_n \in I$ and a corresponding sequence of positive real numbers $t_1, \ldots, t_n$ such that
$$ 
\Phi_{i_n}^{t_n} \circ \ldots \circ \Phi_{i_1}^{t_1}(y) = x; 
$$ 
\item \emph{accessible from $y$} if for any neighborhood $U$ of $x$ there is $z \in U$ such that $z$ is reachable from $y$;  
\item \emph{accessible from $S \subset \R^d$} if it is accessible from any $y \in S$. If $x$ is accessible from $\R^d$, we simply say that $x$ is accessible. 
\end{enumerate} 
\end{defn}

A point $x \in \R^d$ is accessible if and only if for every neighborhood $U$ of $x$, for every $y \in \R^d$ and for every $i, j \in I$ there is $t > 0$ such that $\Pp^t_{y,i}(U \times \{j\}) > 0$.  If $x$ is accessible, then the points $(x,i)$, $i \in I$, are contained in the support of any IPM for $(\Pp^t)$. 

\begin{thm}[Bena\"im, Strickler,~\cite{BS17}]     \label{thm:persistence}
Let $M_+ = M \setminus \{0\}$. The following statements hold. 
\begin{enumerate}
\item If $\Lambda^{-} > 0$, then there exists an IPM $\mu$ of $(X,E)$ such that $\mu(M_+ \times I) = 1$. In addition, for any starting point $x \in M_+$, $X_t$ almost surely does not converge to $0$ as $t \to \infty$. 
\item If $\Lambda^+ < 0$ and if the point $0$ is accessible, then for any starting point $x \in M$, $X_t$ converges almost surely to $0$ as $t \to \infty$. In particular, there is no IPM that assigns positive mass to $M_+ \times I$. 
\end{enumerate}
\end{thm}

Now, we review sufficient conditions for uniqueness and absolute continuity of the IPM. Recall that the \emph{Lie bracket} of $C^1$ vector fields $f_0$ and $f_1$ on $\R^d$ is defined as 
$$ 
[f_0,f_1](x) = Df_1(x) f_0(x) - Df_0(x) f_1(x), \quad x \in \R^d. 
$$ 
Let $\cL$ denote the Lie algebra generated by $(f_i)_{i \in I}$, i.e. $\cL$ is the smallest collection of $C^{\infty}$ vector fields on $\R^d$ that contains $(f_i)_{i \in I}$, and is closed under linear combinations and the Lie bracket operation. 

\begin{defn}      \label{def:Hoermander}
We say that the \emph{weak bracket condition} is satisfied at a point $x \in \R^d$ if 
$$ 
\{f(x): f \in \cL\} = \R^d. 
$$ 
\end{defn}

The weak bracket condition is essentially H\"ormander's condition for smoothness of transition densities for a diffusion process with the noise acting along $(f_i)_{i \in I}$, see~\cite[Section 2.3]{Nualart}.  

\begin{thm}[Bena\"im, Le Borgne, Malrieu, Zitt,~\cite{Benaim}; Bakhtin, Hurth,~\cite{Bakhtin}]   \label{thm:unique_im}
Let $U \subset \R^d$ be an open positive invariant set. Suppose $(\Pp^t)$ admits an IPM $\mu$ such that $\mu(U \times I) = 1$. Assume in addition that there exists $x \in U$ such that (i) $x$ is accessible from $U$ and (ii) the weak bracket condition holds at $x$.  Then, $\mu$ is the unique IPM assigning full measure to $U \times I$, and $\mu$ is absolutely continuous with respect to the product of Lebesgue measure on $\R^d$ and counting measure on $I$.  
\end{thm}

If $d=1$, the weak bracket condition holds at any point that is not an equilibrium of all $(f_i)_{i \in I}$. The interesting condition is then existence of an accessible point. 

Suppose now that $(\Pp^t)$ admits an absolutely continuous IPM with probability density function $\rho(x,i)$.  We refer to the projections $\rho_i = \rho(\cdot, i), i \in I,$ as \emph{invariant densities}.  For some simple PDMPs on $\R \times I$, we can give explicit formulas for invariant densities. Besides, we have the following regularity result. 

\begin{thm}[Bakhtin, Hurth, Mattingly,~\cite{Mattingly}]     \label{thm:smooth_id}
Assume that $(f_i)_{i \in I}$ are $C^{\infty}$ vector fields on $\R$ with locally finite sets of critical points each. Let $x \in \R$ such that $f_i(x) \neq 0$ for every $i \in I$. Then, the invariant densities $(\rho_i)_{i \in I}$ of an absolutely continuous IPM are $C^{\infty}$ smooth at $x$. 
\end{thm}

%%%%%%%%%%%%%%%%%%%%%%%%%%%%%%%%%%%%%%%%%%%%%%%%%%%%%%%%%%%%%%%%%%%%%%%%%%%%%%%%
\section{Two Nontrivial Trapping Regions}
\label{sec:stable}

Given a vector field $f$ on $\R^d$ with flow function $\Phi$ and a set $\cV \subset \R^d$, we call $\cV$ a \emph{trapping region} with respect to $f$ if for any $x \in \cV$ and any $t > 0$ we have $\Phi^t(x) \in \cV$.  We split our analysis into two cases, which can occur for our normal forms in
different parameter regimes. Either, there exists a trapping region $\cV\subset \R^d$ 
of finite positive Lebesgue measure. Or, trajectories 
leave any bounded set except for a set of measure zero, which is going to consist of 
unstable equilibria in our case. In this section, we cover the case when such a trapping region exists both below and above 
the bifurcation value. The case when a trapping region exists only below 
or only above the bifurcation value is covered in Section~\ref{sec:unstable}.

%%%%%%%%%%%%%%%%%%%%%%%%%%%%%%%%%%%%%%%%%%%%%%%%
\subsection{Supercritical Pitchfork Bifurcation}
\label{ssec:superpitch}

Consider the ODE~\eqref{eq:ODE} for $d=1$ and assume the existence of a 
trivial branch of equilibria $f(x_*,p)=0$ for all $p$. Assume that the following 
conditions hold at $(x,p)=(x_*,p_*)$:
\be
\partial_x f(x_*,p_*)=0,\quad \partial_{xx}f(x_*,p_*)=0,\quad 
\partial_{xxx}f(x_*,p_*)<0,\quad \partial_{xp}f(x_*,p_*)\neq 0.
\ee 
Then a bifurcation occurs at $(x_*,p_*)$, which can be proven to be locally 
topologically equivalent to the supercritical pitchfork bifurcation
normal form
\be
\label{eq:superpitchforknf}
x'=px-x^3.
\ee
The dynamics of~\eqref{eq:superpitchforknf} is easy to analyze. For $p<0$, there 
is a unique globally stable equilibrium point $x_*=0$. For $p>0$, $x_*=0$ is 
unstable while the equilibria $x_\pm=\pm\sqrt{p}$ are locally stable. For any $p\in\R$, all 
trajectories remain bounded so trapping regions of positive measure are easy to find.
We now analyze the normal form~\eqref{eq:superpitchforknf} from the viewpoint
of PDMPs by switching the parameter $p$.  For fixed parameters $p_{-} < 0$ and $p_{+} > 0$, we switch between the vector fields 
$$ 
f_{-1}(x) = p_{-} x - x^3, \quad f_1(x) = p_+ x - x^3. 
$$ 
We denote the rate of switching from $f_{-1}$ to $f_1$ by $\lambda_{-}$ and the rate of switching from $f_1$ to $f_{-1}$ by $\lambda_+$.  Since $0$ is an equilibrium for both vector fields, the semigroup $(\Pp^t)$ associated with the PDMP $(X,E)$ admits at least one IPM, namely the product of the Dirac measure at $0$ and the measure on $I = \{-1,1\}$ that assigns probability $\tfrac{\lambda_+}{\lambda_+ + \lambda_{-}}$ to $-1$ and $\tfrac{\lambda_{-}}{\lambda_+ + \lambda_{-}}$ to $1$.  The latter is precisely the IPM of the continuous-time Markov chain $E$ on the state space $I$.  For ease of reference, we call this trivial IPM $\delta$. 

\begin{thm}       \label{thm:existence_pitchfork}
The following statements hold. 
\begin{enumerate}
\item If $\tfrac{\lambda_+}{p_+} < -\tfrac{\lambda_{-}}{p_{-}}$, the semigroup $(\Pp^t)$ admits exactly three ergodic IPM: the trivial measure $\delta$, a measure $\mu$ such that $\mu((0,\infty) \times I) = 1$, and a measure $\pi$ such that $\pi((-\infty,0) \times I) = 1$. 
\item If $\tfrac{\lambda_+}{p_+} \geq -\tfrac{\lambda_{-}}{p_{-}}$, then $\delta$ is the unique IPM for $(\Pp^t)$. 
\end{enumerate}
\end{thm}

\begin{thm}     \label{thm:densities_pitchfork}
Suppose that $\tfrac{\lambda_+}{p_+} < -\tfrac{\lambda_{-}}{p_{-}}$.  Then, the ergodic IPM $\mu$ and $\pi$ assigning measure $1$ to $(0,\infty) \times I$ and $(-\infty,0) \times I$, respectively, are absolutely continuous. Moreover, the corresponding invariant densities $\rho^{\mu}$ and $\rho^{\pi}$ are given by   
$$
\rho_i^{\mu}(x) = \rho_i^{\pi}(-x) = C x^{-\frac{\lambda_{-}}{p_{-}} - \frac{\lambda_+}{p_+} - 1} \left(-p_{-} + x^2 \right)^{\frac{\lambda_{-}}{2 p_{-}}  - \frac{1}{2} (1-i)} \left(p_+ - x^2 \right)^{\frac{\lambda_+}{2 p_+} - \frac{1}{2} (1+i)} \mathbbm{1}_{(0,\sqrt{p_+})}(x), \quad i \in I.
$$
Here, $C$ is a normalizing constant. 
\end{thm}

\bpf[Proof of Theorem~\ref{thm:existence_pitchfork}]  Let $M = [0, \sqrt{p_+}]$ and $M_+ = (0,\sqrt{p_+}]$.  Then, $M$ is a compact positive invariant set containing the common equilibrium $0$. Moreover, as $0$ is globally asymptotically stable for $f_{-1}$, $0$ is accessible from $M$.  If we linearize $f_{-1}$ and $f_1$ at $x=0$, we obtain 
$$ 
A_{-1} = \frac{d}{dx} f_{-1}(x) \vert_{x=0} = p_{-}, \quad A_1 = \frac{d}{dx} f_1(x) \vert_{x=0} = p_+. 
$$ 
We have $C_M = [0,\infty)$ and $C_M \cap S^0 = \{1\}$.  The angular process $(\Theta, E)$ has a unique IPM $\nu$ that assigns probability $\tfrac{\lambda_+}{\lambda_+ + \lambda_{-}}$ to $\{1\} \times \{-1\}$ and probability $\tfrac{\lambda_{-}}{\lambda_+ +\lambda_{-}}$ to $\{1\} \times \{1\}$. Thus, 
$$ 
\Lambda^+ = \Lambda^{-} = \Lambda(\nu) = \frac{p_{-} \lambda_+ + p_+ \lambda_{-}}{\lambda_{+} + \lambda_{-}}, 
$$ 
which is positive if $\tfrac{\lambda_+}{p_+} < - \tfrac{\lambda_{-}}{p_{-}}$ and negative if $\tfrac{\lambda_+}{p_+} > -\tfrac{\lambda_{-}}{p_{-}}$.  By Theorem~\ref{thm:persistence}, if $\tfrac{\lambda_+}{p_+} < -\tfrac{\lambda_{-}}{p_{-}}$, there exists an IPM $\mu$ such that $\mu(M_+ \times I) = 1$; and if $\tfrac{\lambda_+}{p_+} > -\tfrac{\lambda_{-}}{p_{-}}$, there is no IPM assigning positive mass to $M_+ \times I$. Suppose now that $\tfrac{\lambda_+}{p_+} < -\tfrac{\lambda_{-}}{p_{-}}$, and consider the open positive invariant set $(0,\infty)$.  Let $x \in (0,\sqrt{p_+})$ and observe that $x$ is accessible from $(0,\infty)$.  Since $f_{-1}$ and $f_1$ do not vanish at $x$, the weak bracket condition is satisfied.  By Theorem~\ref{thm:unique_im}, there is exactly one IPM $\mu$ assigning full measure to $(0,\infty) \times I$. This measure is ergodic: By the ergodic decomposition theorem (see, e.g.,~\cite[Theorem 5.7]{Hairer:ergodicity-lectures}), there exists an ergodic IPM $\pi$ assigning positive mass to $(0,\infty) \times I$. Since $(0,\infty)$ is positive invariant, we have $\pi((0,\infty) \times I) = 1$ and hence $\pi = \mu$. This argument also shows that $\mu$ is the only ergodic IPM that assigns positive mass to $(0,\infty) \times I$.  A completely analogous reasoning applies to the positive invariant set $(-\infty,0)$.  

It remains to consider the critical case $\tfrac{\lambda_+}{p_+} = -\tfrac{\lambda_{-}}{p_{-}}$, where Theorem~\ref{thm:persistence} does not apply. To obtain a contradiction, we assume that there is an ergodic IPM $\mu$ that, without loss of generality, assigns measure $1$ to $(0,\infty) \times I$. By Theorem~\ref{thm:unique_im}, $\mu$ has a density $\rho$, and by Theorem~\ref{thm:smooth_id} $\rho_{-1}$ and $\rho_1$ are smooth in $(0,\sqrt{p_+})$. Therefore, they satisfy the formula in Theorem~\ref{thm:densities_pitchfork}.  As $\tfrac{\lambda_{-}}{p_{-}} + \tfrac{\lambda_+}{p_+} = 0$, $\rho_{-1}$ and $\rho_1$ behave asymptotically as $x^{-1}$ as $x \downarrow 0$. Since $x^{-1}$ is not integrable in a neighborhood of $0$, we arrive at a contradiction. 
\epf 

\bigskip

\bpf[Proof of Theorem~\ref{thm:densities_pitchfork}] Absolute continuity of $\mu$ and $\pi$ follows from Theorem~\ref{thm:unique_im}.  As shown in the proof of Theorem~\ref{thm:existence_pitchfork}, $\mu((0,\sqrt{p_+}] \times I) = 1$, so the invariant densities $(\rho_i^{\mu})_{i \in I}$ vanish outside of $[0,\sqrt{p_+}]$. By Theorem~\ref{thm:smooth_id}, $(\rho^{\mu}_i)_{i \in I}$ are $C^{\infty}$ on $(0,\sqrt{p_+})$ and thus satisfy the Fokker--Planck equations, see for instance~\cite{Faggionato}.  Written in terms of probability fluxes $\varphi_i = \rho^{\mu}_i f_i, i \in I$, the Fokker--Planck equations read for $x \in (0,\sqrt{p_+})$ 
\begin{align}
\varphi_{-1}'(x) =& -\frac{\lambda_{-}}{f_{-1}(x)} \varphi_{-1}(x) + \frac{\lambda_+}{f_1(x)} \varphi_1(x), \label{eq:Fokker_1} \\
\varphi_1'(x) =& -\frac{\lambda_+}{f_1(x)} \varphi_1(x) + \frac{\lambda_{-}}{f_{-1}(x)} \varphi_{-1}(x). \label{eq:Fokker_2}
\end{align}
Then, 
$$ 
\varphi_{-1}' + \varphi_1' \equiv 0, 
$$ 
so $\varphi_{-1} + \varphi_1$ is constant. We even have $\varphi_{-1} + \varphi_1 \equiv 0$~\cite{Balazs}.  The ODE in~\eqref{eq:Fokker_1} becomes 
$$ 
\varphi_{-1}'(x) = -\left(\frac{\lambda_{-}}{f_{-1}(x)} + \frac{\lambda_+}{f_1(x)} \right) \varphi_{-1}(x),  
$$ 
which is solved by 
\begin{align*} 
\varphi_{-1}(x) =& C \exp \biggl( -\lambda_{-} \int \frac{dx}{f_{-1}(x)} - \lambda_+ \int \frac{dx}{f_1(x)} \biggr) \\ 
=& C x^{-\frac{\lambda_{-}}{p_{-}} - \frac{\lambda_+}{p_+}} \left(-p_{-}+x^2 \right)^{\frac{\lambda_{-}}{2p_{-}}} \left(p_+-x^2\right)^{\frac{\lambda_+}{2p_+}}.  
\end{align*}
We obtain the desired formula for $\rho^{\mu}$ with $\rho^{\mu}_{-1} = \varphi_{-1} / f_{-1}$ and $\rho^{\mu}_1 = -\varphi_{-1} / f_1$.  The formula for $\rho^{\pi}$ follows from the fact that both $f_{-1}$ and $f_1$ are odd. 
\epf 

\bigskip

%%%%%%%%%%%%%%%%%%%%%%%%%%%%%%%%%%%%%%%%%%%%%%%%
\subsection{Supercritical Hopf Bifurcation}
\label{ssec:superHopf}

Consider the ODE~\eqref{eq:ODE} for $d=2$. Assume that $x_*=x_*(p)$
is a family of equilibrium points for all $p$ in a parameter-space 
neighbourhood of $p_*$. Let $A=\txtD_xf (x_*(p),p)$ and assume that 
\be
\textnormal{spec}(A)=\{\alpha(p)\pm \txti\omega(p)\},\quad \alpha(p_*)=0,
\quad \alpha'(p_*)\neq 0,\quad \omega(p_*)\neq 0.
\ee 
Furthermore, consider the first Lyapunov coefficient $l_1=l_1(p)$, which
is computable from $f$ using partial derivatives up to and including third
order; see the formulas in~\cite{GH,Kuznetsov}. Assume that $l_1(p_*)<0$.
Then a bifurcation occurs at $(x_*,p_*)$, which can be proven to be locally 
topologically equivalent to the supercritical Hopf bifurcation
normal form
\be
\label{eq:superHopfnf}
\begin{array}{lcl}
x_1'&=&px_1-x_2-x_1(x_1^2+x_2^2),\\
x_2'&=&x_1+px_2-x_2(x_1^2+x_2^2).
\end{array}
\ee
The dynamics of~\eqref{eq:superHopfnf} can be analyzed a lot easier upon 
changing to polar coordinates $(x_1,x_2)=(r\cos \theta,r\sin\theta)$, which
gives 
\be
\label{eq:superHopfpolar}
\begin{array}{lcl}
\theta' &=& 1, \\
r'&=&pr-r^3.
\end{array}
\ee
Analyzing the simple vector field~\eqref{eq:superHopfpolar} and returning 
to Euclidean coordinates, one finds that for $p<0$, there 
is a unique globally stable equilibrium point $x_*=0$. For $p>0$, $x_*=0$ is 
unstable while there exists a family of stable periodic orbits 
$\{\|x\|_2=\sqrt{p}\}$. For any $p\in\R$, all 
trajectories remain bounded so trapping regions of positive measure always exist.
We now analyze the normal form~\eqref{eq:superHopfnf} from the viewpoint
of PDMPs by switching the parameter $p$, again working in polar coordinates.  For fixed $p_{-} < 0$ and $p_+ > 0$, we switch between the vector fields 
$$ 
g_{-1}(\theta, r) = (1, p_{-} r - r^3)^{\top}, \quad g_1(\theta, r) = (1, p_+ r - r^3)^{\top}. 
$$ 
In analogy to the case of the supercritical pitchfork bifurcation, we denote the rate of switching from $g_{-1}$ to $g_1$ by $\lambda_{-}$, and the rate of switching from $g_1$ to $g_{-1}$ by $\lambda_+$.  As before, the origin is an equilibrium for both vector fields, so $\delta$, defined as the product of the Dirac measure at the origin and the discrete measure assigning probability $\tfrac{\lambda_+}{\lambda_+ + \lambda_{-}}$ to $-1$ and $\tfrac{\lambda_{-}}{\lambda_+ + \lambda_{-}}$ to $1$, is an IPM.  Let $\nu$ denote the unique IPM for the PDMP induced by switching between the one-dimensional vector fields 
$$ 
f_{-1}(r) = p_{-} r - r^3, \quad f_1(r) = p_+ r - r^3, \quad r > 0
$$  
at rates $\lambda_{-}$ and $\lambda_+$, whose existence is guaranteed by Theorem~\ref{thm:existence_pitchfork}. 

\begin{thm}     
\label{thm:existence_hopf}
The following statements hold. 
\begin{enumerate}
\item If $\tfrac{\lambda_+}{p_+} < -\tfrac{\lambda_{-}}{p_{-}}$, $(\Pp^t)$ admits exactly two ergodic IPM: the measure $\delta$ and a measure $\mu$ that is the product of  Lebesgue measure on the unit circle $S^1$, normalized by the factor $\tfrac{1}{2 \pi}$, and the IPM $\nu$. 
\item If $\tfrac{\lambda_+}{p_+} \geq -\tfrac{\lambda_{-}}{p_{-}}$, then $\delta$ is the unique IPM for $(\Pp^t)$. 
\end{enumerate}
\end{thm}

\bpf Suppose first that $\tfrac{\lambda_+}{p_+} < -\tfrac{\lambda_{-}}{p_{-}}$.  Let $\mu$ denote the product of Lebesgue measure on $S^1$, normalized by the factor $\tfrac{1}{2 \pi}$, and the IPM $\nu$. Then, for $t > 0$, $i, j \in I$, $\theta \in S^1$, $r > 0$ and measurable sets $A \subset S^1$, $B \subset (0,\infty)$, we have  
$$ 
\Pp^t_{\theta, r, j}(A \times B \times \{i\}) = \mathbbm{1}_A(\theta + t) \widehat \Pp^t_{r,j}(B \times \{i\}), 
$$ 
where $\theta + t$ should be understood modulo $2 \pi$, and where $\widehat \Pp$ denotes the semigroup associated with the PDMP induced by $f_{-1}$ and $f_1$.  This form of independence for $\theta$ and $r$ holds because the evolution of the angular component $\theta$ is entirely deterministic and in particular not affected by the switching times.  Thus, 
\begin{align*}
\mu \Pp^t(A \times B \times \{i\}) =& \frac{1}{2 \pi} \int_{S^1} \mathbbm{1}_A(\theta + t) \ d \theta \sum_{j \in \{-1,1\}} \int_0^{\infty} \widehat \Pp^t_{r, j}(B \times \{i\}) \ \nu(dr \times \{j\})   \\
=& \frac{\Leb(A)}{2 \pi} \nu \widehat \Pp^t(B \times \{i\}) = \frac{\Leb(A)}{2 \pi} \nu (B \times \{i\}) = \mu(A \times B \times \{i\}).  
\end{align*}
Hence, $\mu$ is an IPM for $(\Pp^t)$. Next, we show that $\mu$ is the only IPM such that $\mu(S^1 \times (0,\infty) \times I) = 1$.  First we show that any point in $S^1 \times (0, \sqrt{p_+})$ is accessible from $S^1 \times (0,\infty)$. Fix two points $(\alpha, p) \in S^1 \times (0,\sqrt{p_+})$ and $(\beta, q) \in S^1 \times (0,\infty)$. For $i \in I$, we denote the flow associated with the vector field $g_i$ by $\Phi_i$.  As $s \to \infty$, the radial component of $\Phi_{-1}^s(\beta,q)$ tends to $0$.  Let $s > 0$ such that $\Phi_{-1}^s(\beta,q)$ has angular component $\alpha$ and radial component $u < p$. As $p \in (0,\sqrt{p_+})$, a short computation shows that $\Phi_1^t(\alpha, p)$ has radial component  
$$ 
\left(\frac{e^{2 p_+ t} p_+ C}{1 + e^{2 p_+ t} C} \right)^{\frac{1}{2}}, 
$$ 
where $C = \tfrac{p_+ p^2}{p_+^2 - p_+ p^2} > 0$. This shows that the vector field $g_1$ is both forward and backward complete on the punctured disk $S^1 \times (0,\sqrt{p_+})$, with $\lim_{t \to -\infty} \Phi_1^t(\alpha, p) = 0$.  Hence, there is $t < 0$ such that $\Phi_1^t(\alpha, p)$ has angular component $\alpha$ and radial component $v < u$. For $T \geq 0$, let $h(T)$ denote the difference of the radial components of $\Phi_{-1}^T(\alpha, u)$ and $\Phi_1^T(\alpha,v)$. Then, $h(0) = u - v > 0$ and $h(-t) \leq u - p < 0$. As $h$ is continuous, there is $\tau \in (0,-t)$ such that $h(\tau) = 0$.  As the points $\Phi_{-1}^T(\alpha,u)$ and $\Phi_1^T(\alpha, v)$ have the same angular component for any $T \geq 0$, we have 
$$ 
\Phi_{-1}^{\tau}(\alpha, u) = \Phi_1^{\tau}(\alpha, v). 
$$ 
Thus, we can reach the point $(\alpha, p)$ from $(\beta, q)$ as follows: First, flow along the vector field $g_{-1}$ for time $s + \tau$, then make a switch and flow along $g_1$ for time $-t - \tau$. 

For $(\alpha, p) \in S^1 \times (0, \sqrt{p_+})$, the vectors $g_{-1}(\alpha, p)$ and $g_1(\alpha, p)$ are clearly transversal, so the weak bracket condition holds as well.  By Theorem~\ref{thm:unique_im}, $\mu$ is indeed the only IPM assigning mass $1$ to $S^1 \times (0, \infty)$. The fact that $\delta$ and $\mu$ are the only ergodic IPM follows along the same lines as in the proof of Theorem~\ref{thm:existence_pitchfork}. 

Now, we consider the case $\tfrac{\lambda_+}{p_+} \geq -\tfrac{\lambda_{-}}{p_{-}}$. To obtain a contradiction, suppose that there is an IPM $\pi$ for $(\Pp^t)$ such that $\pi(S^1 \times (0,\infty) \times I) > 0$. By the ergodic decomposition theorem, we may assume without loss of generality that $\pi(S^1 \times (0,\infty) \times I) = 1$. Consider the marginal $\widehat \pi(\cdot) = \pi(S^1 \times \cdot)$, which is a probability measure on $(0,\infty) \times I$. For $t > 0$ and with $\widehat \Pp$ defined as above, we have for measurable $B \subset (0,\infty)$ and $i \in I$ 
\begin{align*} 
\widehat \pi \widehat \Pp^t(B \times \{i\}) =& \sum_{j \in I} \int_0^{\infty} \widehat \Pp^t_{r,j}(B \times \{i\}) \ \widehat \pi(dr \times \{j\}) \\
=& \sum_{j \in I} \int_{S^1} \int_0^{\infty} \widehat \Pp^t_{r,j}(B \times \{i\}) \  \pi(d \theta \times dr \times \{j\}) \\
=& \sum_{j \in I} \int_{S^1} \int_0^{\infty} \mathbbm{1}_{S^1}(\theta + t)  \widehat \Pp^t_{r,j}(B \times \{i\}) \ \pi(d \theta \times dr \times \{j\}) \\
=& \sum_{j \in I} \int_{S^1}\int_0^{\infty} \Pp^t_{\theta, r, j}(S^1 \times B \times \{i\}) \ \pi(d \theta \times dr \times \{j\}) = \pi(S^1 \times B \times \{i\}) = \widehat \pi(B \times \{i\}). 
\end{align*}  
This computation shows that $\widehat \pi$ is an IPM for $(\widehat \Pp^t)$. But Theorem~\ref{thm:existence_pitchfork} implies that $(\widehat \Pp^t)$ has no IPM if $\tfrac{\lambda_+}{p_+} \geq -\tfrac{\lambda_{-}}{p_{-}}$, a contradiction. 
\epf 

\bigskip  

%%%%%%%%%%%%%%%%%%%%%%%%%%%%%%%%%%%%%%%%%%%%%%%%
\subsection{Transcritical Bifurcation}
\label{ssec:transcritical}

Consider the ODE~\eqref{eq:ODE} for $d=1$ and assume the existence of a 
trivial branch of equilibria $f(x_*,p)=0$ for all $p$. Assume that the following 
conditions hold at $(x,p)=(x_*,p_*)$:
\be
\partial_x f(x_*,p_*)=0,\quad \partial_{xx}f(x_*,p_*)\neq 0,\quad 
\partial_{xp}f(x_*,p_*)\neq 0.
\ee 
Then a bifurcation occurs at $(x_*,p_*)$, which can be proven to be locally 
topologically equivalent to the transcritical bifurcation
normal form
\be
\label{eq:transcriticalnf}
x'=px-x^2.
\ee
The dynamics of~\eqref{eq:transcriticalnf} works as follows. There 
are two families of equilibrium points $x_*=0$ and $x_{**}=p$. For $p<0$, 
$x_*$ is locally stable, while $x_{**}$ is unstable. For $p>0$, the stabilities 
switch. There are bounded trapping regions of positive measure given in the different 
parameter regimes by
\benn
\cV_{p<0}=[x_{**},0]\qquad\text{and}\qquad \cV_{p>0}=[0,x_{**}] 
\eenn
with the special case $\cV_{p=0}=[0,K]$ for any $K>0$. For fixed $p_{-} < 0$ and $p_+ > 0$, consider the vector fields
$$ 
f_{-1}(x) = p_{-} x - x^2, \quad f_1(x) = p_+ x - x^2. 
$$ 
These vector fields are not forward complete: trajectories for $f_{-1}$ that start to the left of $p_{-}$ and trajectories for $f_1$ that start to the left of $0$ move off to $-\infty$ in finite time.  To obtain a well-defined PDMP, we therefore restrict ourselves to the positive invariant set $[0,\infty)$, where both $f_{-1}$ and $f_1$ have bounded trajectories and are in particular forward complete. We switch from $f_{-1}$ to $f_1$ at rate $\lambda_{-}$ and from $f_1$ to $f_{-1}$ at rate $\lambda_+$, and we let $\delta$ denote the product of the Dirac measure at $0$ and the IPM of $E$. 

\begin{remark}
It is possible to define a PDMP that involves switching between $f_{-1}$ and $f_1$ on the larger interval $(p_{-},\infty)$.  Since $(p_{-}, \infty)$ is not a trapping region for $f_1$, one needs to ensure that we switch away from $f_1$ before reaching the point $p_{-}$. This can be achieved by letting the switching rate $\lambda_+$ depend on the location $x$ of the switching trajectory, with $\lambda_+(x)$ blowing up as $x \to p_{-}$ from the right. 
\end{remark}

\begin{thm}   \label{thm:existence_transcritical}
The following statements hold. 
\begin{enumerate}
\item If $\tfrac{\lambda_+}{p_+} < -\tfrac{\lambda_{-}}{p_{-}}$, there are exactly two ergodic IPM: the trivial measure $\delta$ and a measure $\mu$ such that $\mu((0,\infty) \times I) = 1$. 
\item If $\tfrac{\lambda_+}{p_+} \geq -\tfrac{\lambda_{-}}{p_{-}}$, then $\delta$ is the only IPM. 
\end{enumerate}
\end{thm} 

This statement can be shown along the same lines as Theorem~\ref{thm:existence_pitchfork}. We therefore omit the proof.

\begin{thm}     \label{thm:densities_transcritical}
Suppose that $\tfrac{\lambda_+}{p_+} < -\tfrac{\lambda_{-}}{p_{-}}$. Then, the ergodic IPM $\mu$ assigning mass $1$ to $(0,\infty) \times I$ is absolutely continuous. Moreover, the corresponding invariant density $\rho$ is given by   
$$ 
\rho_i(x) = C x^{-\frac{\lambda_{-}}{p_{-}} - \frac{\lambda_+}{p_+}-1} (-p_{-}+x)^{\frac{\lambda_{-}}{p_{-}}- \frac{1}{2}(1-i)} (p_+ - x)^{\frac{\lambda_+}{p_+}-\frac{1}{2}(1+i)} \mathbbm{1}_{(0,p_+)}(x), \quad i \in I.  
$$
\end{thm}

\bpf Absolute continuity of $\mu$ follows from Theorem~\ref{thm:unique_im}. As $\mu((0,p_+] \times I) = 1$, the invariant densities $(\rho_i)_{i \in I}$ vanish outside of $[0,p_+]$. By Theorem~\ref{thm:smooth_id}, $(\rho_i)_{i \in I}$ are $C^{\infty}$ on $(0,p_+)$ and thus satisfy the Fokker -- Planck equations. For the probability flux $\varphi_{-1}$, we have 
$$
\varphi_{-1}(x) = C \exp\left(- \lambda_{-} \int \frac{dx}{f_{-1}(x)} - \lambda_+ \int \frac{dx}{f_1(x)} \right) =  C x^{-\frac{\lambda_{-}}{p_{-}} - \frac{\lambda_+}{p_+}}  \left(- p_{-}+x\right)^{\frac{\lambda_{-}}{p_{-}}} \left(p_+ - x\right)^{\frac{\lambda_+}{p_+}}.  
$$ 
As in the case of the supercritical pitchfork bifurcation, we obtain the desired formula with $\rho_{-1} = \varphi_{-1} / f_{-1}$ and $\rho_1 = - \varphi_{-1} / f_1$. 
\epf

\medskip

If the switching rates $\lambda_+$ and $\lambda_{-}$ do not depend on $X$, the PDMP $(X,E)$ starting at a point to the left of $0$ will tend to $-\infty$ in finite time with positive probability.  To make this statement more precise, we define for $a < p_{-}$ the stopping time
$$ 
\tau_a = \inf \{t \geq 0: X_t \leq a\}. 
$$

\begin{prop}     \label{prop:divergence_to_infinity}
Let $\nu$ be a probability measure on $\R \times I$ such that $\nu((-\infty,0) \times I) = 1$, and let $\nu \Pp^{\tau_a}$ denote the law of $(X,E)$ with initial distribution $\nu$ and stopped at time $\tau_a$.  There is a nonincreasing function $g: (-\infty, p_{-}-2] \to (0,\infty)$ such that $\int_{-\infty}^{p_{-}-2} g(a) \ da < \infty$ and   
$$ 
\nu \Pp^{\tau_a} \left(\tau_{p_{-}-1} < \infty \right) > 0, \quad \nu \Pp^{\tau_a} \left(\tau_a - \tau_{a+1} < g(a) \mid \tau_{p_{-}-1} < \infty\right) = 1, \quad a \leq p_{-} - 2. 
$$ 
\end{prop}

Proposition~\ref{prop:divergence_to_infinity} essentially says that $X_t$ goes off to $-\infty$ in finite time with positive probability if the initial distribution assigns full measure to $(-\infty,0) \times I$:  There is a positive probability that $X$ reaches the interval $(-\infty, p_{-} - 1]$ in finite time. And once $X$ is in $(-\infty, p_{-}-1]$, it blows up to $-\infty$ with probability $1$ in time less than 
$$ 
\left(\tau_{p_{-}-2} - \tau_{p_{-}-1} \right) + \left(\tau_{p_{-}-3} - \tau_{p_{-}-2} \right) + \ldots \leq  \sum_{k =2}^{\infty} g(p_{-}-k) < \infty. 
$$ 

\bpf[Proof of Proposition~\ref{prop:divergence_to_infinity}] Fix $a \leq p_{-}-2$.  Let us first show that $\nu \Pp^{\tau_a} (\tau_{p_{-}-1} < \infty) > 0$.
Let $\delta > 0$ be so small that $\nu((-\infty, -\delta] \times I) > 0$, and let $r, s > 0$ such that 
$$ 
\Phi_{-1}^r(-\delta) = -\frac{\delta}{2}, \quad \Phi_1^s(-\tfrac{\delta}{2}) = p_{-} - 1. 
$$
If $\nu \Pp^{\tau_a}(s + r > \tau_a) > 0$, we also have
$$ 
\nu \Pp^{\tau_a}(\tau_{p_{-}-1} < \infty) > 0. 
$$ 
If $\nu \Pp^{\tau_a}(s+r \leq \tau_a) = 1$, we use the estimate
\begin{equation}     \label{eq:lower_bd_transcritical}
\nu \Pp^{\tau_a}(\tau_{p_{-}-1} < \infty) \geq \nu \Pp^{\tau_a}(\tau_{p_{-}-1} < \infty, X_0 \leq -\delta,  E_t = 1 \ \forall t \in [r, r+s]). 
\end{equation}
Suppose that $s +r \leq \tau_a$ and $X_0 \leq -\delta$. Then, we have $X_r \leq -\tfrac{\delta}{2}$. If in addition $E_t = 1$ for all $t \in [r,r+s]$, it follows that $\tau_{p_{-}-1} < \infty$. Hence, the term on the right side of~\eqref{eq:lower_bd_transcritical} equals 
$$
\nu \Pp^{\tau_a}(X_0 \leq -\delta, E_t = 1 \ \forall t \in [r,r+s]) > 0.
$$  
Now, we come to the second statement. We will specify the function $g$ later in the proof. Since there is $c > 0$ such that $f_1(x) \leq f_{-1}(x) \leq -c$ for all $x \in (-\infty, p_{-}-1]$, we have $\tau_a < \infty$ for all $a \leq p_{-}-2$ whenever $\tau_{p_{-}-1} < \infty$.  By the strong Markov property, 
$$ 
\nu \Pp^{\tau_a} (\tau_a - \tau_{a+1} < g(a) \mid \tau_{p_{-}-1} < \infty) = \pi \Pp^{\tau_a}(\tau_a < g(a)), 
$$ 
where $\pi$ is the distribution of $(X,E)_{\tau_{a+1}}$ under $\nu \Pp^{\tau_a}(\cdot \mid \tau_{p_{-}-1} < \infty)$, and thus satisfies
$$
\pi((-\infty, a+1] \times I) = 1.
$$ 
In light of $f_1 \leq f_{-1} \leq -c$, we have under $\pi \Pp^{\tau_a}$ 
$$ 
X_t \leq \Phi_{-1}^t(a+1), \quad t \geq 0. 
$$ 
As a result, if we let $g(a)$ be defined by the relation $\Phi_{-1}^{g(a)}(a+1) = a$, we have $\tau_a < g(a)$ under $\pi \Pp^{\tau_a}$.  Since trajectories of $f_{-1}$ starting in $(-\infty, p_{-}-1]$ tend to $-\infty$ in finite time, we also have $\int_{-\infty}^{p_{-}-2} g(a) \ da < \infty$. 
\epf 

\medskip

Proposition~\ref{prop:divergence_to_infinity} and the ergodic decomposition theorem imply that there is no IPM assigning positive mass to $(-\infty,0) \times I$. Looking at Proposition~\ref{prop:divergence_to_infinity}, it is natural to ask under which conditions a blow-up of $X_t$ to $-\infty$ in finite time happens almost surely. The answer follows from Theorem~\ref{thm:divergence_to_infinity} below. 

\begin{thm}      \label{thm:divergence_to_infinity}
Let $a \leq p_{-}-2$ and let $\nu$ be a probability measure on $\R \times I$ such that $\nu((-\infty,0) \times I) = 1$.   
\begin{enumerate}
\item If $\tfrac{\lambda_+}{p_+} < -\tfrac{\lambda_{-}}{p_{-}}$, we have $\nu \Pp^{\tau_a}(\tau_{p_{-}-1} < \infty) = 1$. 

\item If $\tfrac{\lambda_+}{p_+} > -\tfrac{\lambda_{-}}{p_{-}}$ and if $\nu((p_{-},0) \times I) > 0$, we have $\nu \Pp^{\tau_a}(\tau_{p_{-}-1} < \infty) < 1$ and 
$$
\nu \Pp^{\tau_a} \left(\{\tau_{p_{-}-1} < \infty \} \cup \left\{\lim_{t \to \infty} X_t = 0 \right\} \right) = 1.
$$ 
\end{enumerate}
\end{thm}

As stated in the following lemma, with probability $1$, $X_t$ either diverges to $-\infty$ in finite time or converges to $0$ as $t \to \infty$. 

\begin{lem}     \label{lm:divergence_to_infinity}
If $a \leq p_{-}-2$ and if $\nu$ is a probability measure on $\R \times I$ such that $\nu((-\infty,0) \times I) = 1$, we have
$$ 
\nu \Pp^{\tau_a} \left( \{\tau_{p_{-}-1} < \infty\} \cup \left\{\lim_{t \to \infty} X_t = 0 \right\} \right) = 1. 
$$ 
\end{lem}
   
\bpf[Proof of Lemma~\ref{lm:divergence_to_infinity}]  Under $\nu \Pp^{\tau_a}$, the complement of $\{\tau_{p_{-}-1} < \infty\} \cup \{\lim_{t \to \infty} X_t = 0\}$ is 
$$ 
\{\tau_{p_{-}-1} = \infty\} \cap \bigcup_{n=1}^{\infty} \bigcap_{k=1}^{\infty} \bigcup_{t \geq k} \left\{X_t \leq - \frac{1}{n} \right\}. 
$$
For fixed $n \in \N$, consider the event 
$$ 
\bigcap_{k=1}^{\infty} \bigcup_{t \geq k} \left\{X_t \leq -\frac{1}{n} \right\}. 
$$ 
On this event, there is $T > 0$ such that $X_t \leq -\tfrac{1}{n}$ for every $t \geq T$, or there is a sequence of times $t_j \uparrow \infty$ such that $(X,E)_{t_j} = (-\tfrac{1}{n}, 1)$ for every $j$. In the former case, let $s > 0$ such that 
$$ 
\Phi_1^s(-\tfrac{1}{n}) = p_{-} - 1. 
$$
Then, $\tau_{p_{-}-1} < \infty$ or, $\nu \Pp^{\tau_a}$-almost surely, there is $r > T$ such that $E_t = 1$ for $r \leq r \leq r + s$, which also yields $\tau_{p_{-}-1} < \infty$. In the latter case, observe that the first return time to state $(-\tfrac{1}{n},1)$ is finite with probability strictly less than $1$, so by the strong Markov property the event $\{(X,E)_{t_j} = (-\tfrac{1}{n},1) \ \forall j\}$ has probability $0$. 
\epf 

\medskip 

\bpf[Proof of Theorem~\ref{thm:divergence_to_infinity}] Assume first that $\tfrac{\lambda_+}{p_+} < -\tfrac{\lambda_{-}}{p_{-}}$. By Lemma~\ref{lm:divergence_to_infinity} it suffices to show that 
\begin{equation}    \label{eq:no_conv_to_zero} 
\nu \Pp^{\tau_a} \left(\lim_{t \to \infty} X_t = 0 \right) = 0. 
\end{equation} 
Let $M = [\tfrac{p_{-}}{2},0]$ and $M_+ = [\tfrac{p_{-}}{2},0)$. Let $\tilde f_1$ be a smooth vector field that coincides with $f_1$ on the interval $[\tfrac{p_{-}}{4},0]$, is strictly negative on $(\tfrac{p_{-}}{2},0)$, and has $\tfrac{p_{-}}{2}$ as an equilibrium point. In addition, we assume that $\tilde f_1(x) \geq f_1(x)$ for all $x \in \R$.  Then, $M$ is positive invariant for the vector fields $f_{-1}$ and $\tilde f_1$, and $0$ is accessible from $M$. Let $(\tilde X, \tilde E)$ denote the PDMP with vector fields $f_{-1}, \tilde f_1$ and switching rates $\lambda_{-}, \lambda_+$.     Following the proof of Theorem~\ref{thm:existence_pitchfork} and applying Theorem~\ref{thm:persistence}, we see that for any starting point $x \in M_+$, $\tilde X_t$ almost surely does not converge to $0$ as $t \to \infty$.  The Markov property and the fact that any switching trajectory starting in $(-\infty,0)$ and converging to $0$ has to visit points in $M_+$ imply that this result extends to starting points $x \in (-\infty,0)$. Since $\tilde f_1 \geq f_1$, we finally infer~\eqref{eq:no_conv_to_zero}. 

Now, we consider the case $\tfrac{\lambda_+}{p_+} >  -\tfrac{\lambda_{-}}{p_{-}}$, assuming that $\nu((p_{-},0) \times I) > 0$. Defining $\tilde f_1$ and $(\tilde X, \tilde E)$ as above, we obtain with Theorem~\ref{thm:persistence} that for any starting point $x \in M = [\tfrac{p_{-}}{2},0]$, $\tilde X_t$ converges almost surely to $0$ as $t \to \infty$. Fix $x \in [\tfrac{p_{-}}{4},0)$. Then,
$$ 
\delta_{x,-1} \tilde \Pp \left(\lim_{t \to \infty} \tilde X_t = 0 \right) = 1, 
$$ 
where $\delta_{x,-1} \tilde \Pp$ is the distribution of $(\tilde X, \tilde E)$ with initial distribution $\delta_{x,-1}$. The first return time for state $(x,-1)$ must then be infinite with positive probability. In other words, there is a positive probability that the PDMP $(\tilde X, \tilde E)$ starting in $(x,-1)$ stays in $[\tfrac{p_{-}}{4},0) \times I$ for all $t \geq 0$ and thus coincides with $(X,E)$ starting in $(x,-1)$.  In particular, $\Pp_{x,-1}^{\tau_a}(\lim_{t \to \infty} X_t = 0) > 0$ for all $x \in (\tfrac{p_{-}}{4},0)$. Let $\eps > 0$ be so small that $\nu((p_{-}+\eps, 0) \times I) > 0$, and let $s > 0$ such that $\Phi_{-1}^s(p_{-}+\eps) = \tfrac{p_{-}}{4}$. Then, for any $(y,i) \in (p_{-}+\eps,0)$, $(\Pp^s_{y,i})^{\tau_a}((\tfrac{p_{-}}{4},0) \times \{-1\}) > 0$. It follows that 
$$ 
\nu \Pp^{\tau_a} \left(\lim_{t \to \infty} X_t = 0 \right) \geq \sum_{i \in I} \int_{p_{-}+\eps}^0 \int_{\frac{p_{-}}{4}}^0  \Pp^{\tau_a}_{x,-1} \left(\lim_{t \to \infty} X_t = 0 \right) (\Pp_{y,i}^s)^{\tau_a}(dx \times \{-1\}) \ \nu(dy \times \{i\}) > 0.
$$ 
The claim made in part 2 of Theorem~\ref{thm:divergence_to_infinity} then follows from Lemma~\ref{lm:divergence_to_infinity}.  
\epf 

\medskip

%%%%%%%%%%%%%%%%%%%%%%%%%%%%%%%%%%%%%%%%%%%%%%%%%%%%%%%%%%%%%%%%%%%%%%%%%%%%%%%%
\section{One Nontrivial Trapping Region}
\label{sec:unstable}

%%%%%%%%%%%%%%%%%%%%%%%%%%%%%%%%%%%%%%%%%%%%%%%%
\subsection{Subcritical Pitchfork Bifurcation}
\label{ssec:subpitch}

Consider the ODE~\eqref{eq:ODE} for $d=1$ and assume the existence of a 
trivial branch of equilibria $f(x_*,p)=0$ for all $p$. Assume that the following 
conditions hold at $(x,p)=(x_*,p_*)$:
\be
\partial_x f(x_*,p_*)=0,\quad \partial_{xx}f(x_*,p_*)=0,\quad 
\partial_{xxx}f(x_*,p_*)>0,\quad \partial_{xp}f(x_*,p_*)\neq 0.
\ee 
Then a bifurcation occurs at $(x_*,p_*)$, which can be proven to be locally 
topologically equivalent to the subcritical pitchfork bifurcation
normal form
\be
\label{eq:subpitchforknf}
x'=px+x^3.
\ee
The dynamics of~\eqref{eq:subpitchforknf} works as follows. For $p<0$, there 
are three equilibrium points $x_*=0$ and $x_\pm=\pm\sqrt{-p}$. $x_*$ is locally
stable, while $x_\pm$ are unstable. For $p>0$, $x_*=0$ is the only equilibrium
point and it is unstable. For $p\geq 0$ there is no trapping 
region of positive measure. However, $[x_-,x_+]=\cV$ is a trapping region for the 
dynamics when $p<0$.  For fixed $p_{-} < 0$ and $p_+ > 0$, we switch between 
$$ 
f_{-1}(x) = p_{-} x + x^3, \quad f_1(x) = p_+ x + x^3
$$ 
at rates $\lambda_{-}$ and $\lambda_+$. The trivial measure $\delta$ is defined exactly as for the supercritical pitchfork bifurcation. As for the transcritical bifurcation, the vector fields $f_{-1}$ and $f_1$ are not forward complete. E.g., any trajectory of $f_1$ not starting at $0$ blows up in finite time.  If one wishes to define a PDMP outside of the common equilibrium point $0$, one can either let the rate $\lambda_+$ of switching from $f_1$ to $f_{-1}$ depend on the location $x$, with $\lambda_+(x)$ diverging to $\infty$ as $x$ approaches $-\sqrt{-p_{-}}$ from the right and $\sqrt{-p_{-}}$ from the left. Or one can stop the PDMP with constant switching rates once it reaches certain thresholds.  For the latter model, $\delta$ is the unique IPM.  Besides, we have the following result that is reminiscent of Proposition~\ref{prop:divergence_to_infinity} and Theorem~\ref{thm:divergence_to_infinity}.  Since $f_{-1}$ and $f_1$ are odd functions, we may restrict ourselves to the interval $(0,\infty)$, with the understanding that there are completely analogous statements about $(-\infty,0)$. 

\begin{thm}     \label{thm:divergence_subpitch}
Let $\nu$ be a probability measure on $\R \times I$ such that $\nu((0,\infty) \times I) = 1$, and let 
$$ 
\tau_a = \inf\{t \geq 0:  X_t  \geq a\} 
$$ 
for $a > \sqrt{-p_{-}}$. Let $\nu \Pp^{\tau_a}$ denote the law of $(X,E)$ with initial distribution $\nu$ and stopped at time $\tau_a$.  
\begin{enumerate}
\item There is a nonincreasing function $g: [\sqrt{-p_{-}} + 2, \infty) \to (0,\infty)$ such that $\int_{\sqrt{-p_{-}}+2}^{\infty} g(a) \ da < \infty$ and 
$$ 
\nu \Pp^{\tau_a}(\tau_{\sqrt{-p_{-}}+1} < \infty) > 0, \quad \nu \Pp^{\tau_a}(\tau_a - \tau_{a-1} < g(a) \mid \tau_{\sqrt{-p_{-}}+1} < \infty) = 1, \quad a \geq \sqrt{-p_{-}}+2. 
$$ 
\item If $\tfrac{\lambda_+}{p_+} < -\tfrac{\lambda_{-}}{p_{-}}$, we have $\nu \Pp^{\tau_a}(\tau_{\sqrt{-p_{-}}+1} < \infty) = 1$ for $a \geq \sqrt{-p_{-}}+2$. 
\item If $\tfrac{\lambda_+}{p_+} > -\tfrac{\lambda_{-}}{p_{-}}$ and $\nu((0,\sqrt{-p_{-}}) \times I) > 0$, we have $\nu \Pp^{\tau_a}(\tau_{\sqrt{-p_{-}}+1} < \infty) < 1$ and 
$$ 
\nu \Pp^{\tau_a} \left( \{\tau_{\sqrt{-p_{-}}+1} < \infty\} \cup \left\{\lim_{t \to \infty} X_t = 0 \right\} \right) = 1. 
$$  
\end{enumerate}
\end{thm}

The proof is analogous to the ones of Proposition~\ref{prop:divergence_to_infinity} and Theorem~\ref{thm:divergence_to_infinity}, and we omit it.

 \medskip

%%%%%%%%%%%%%%%%%%%%%%%%%%%%%%%%%%%%%%%%%%%%%%%%
\subsection{Subcritical Hopf Bifurcation}
\label{ssec:subHopf}

Consider the same setting as in Section~\ref{ssec:superHopf}, except that
we now assume that the first Lyapunov coefficient satisfies $l_1(p_*)>0$.
This leads to a subcritical Hopf bifurcation normal form
\be
\label{eq:subHopfnf}
\begin{array}{lcl}
x_1'&=&px_1-x_2+x_1(x_1^2+x_2^2),\\
x_2'&=&x_1+px_2+x_2(x_1^2+x_2^2).
\end{array}
\ee
Here the unstable bifurcating family of periodic orbits $\{\|x\|_2=\sqrt{-p}\}$
exists for $p<0$ and in this case $x_*=0$ is locally stable. $x_*$ is 
unstable for $p\geq 0$. For $p<0$, there is a trapping region of positive measure $\cV_{p<0}=\{x\in\R^2:\|x\|_2\leq \sqrt{-p}\}$. 
After a change of variables to polar coordinates, the system in~\eqref{eq:subHopfnf} becomes 
\begin{align*}
\theta' =& 1, \\
r' =& pr + r^3. 
\end{align*}
For fixed $p_{-} < 0$ and $p_+ > 0$, we then switch between 
$$ 
g_{-1}(\theta,r) = (1,p_{-} r+r^3)^{\top}, \quad g_1(\theta,r) = (1,p_+ r+r^3)^{\top}
$$ 
at rates $\lambda_{-}$ and $\lambda_+$.  Here, we encounter the same issue as for the subcritical pitchfork bifurcation. Then, Theorem~\ref{thm:divergence_subpitch} applies to the PDMP induced by the vector fields 
$$ 
f_{-1}(r) = p_{-} r + r^3, \quad f_1(r) = p_+ r + r^3, 
$$ 
and thus to the evolution of the radial component of the PDMP induced by $g_{-1}$ and $g_1$.

\medskip

%%%%%%%%%%%%%%%%%%%%%%%%%%%%%%%%%%%%%%%%%%%%%%%%
\subsection{Fold Bifurcation}
\label{ssec:fold}

Consider the ODE~\eqref{eq:ODE} for $d=1$. Assume that the following 
conditions hold at $(x,p)=(x_*,p_*)$:
\be
f(x_*,p_*)=0=\partial_x f(x_*,p_*)=0,\quad \partial_{xx}f(x_*,p_*)\neq 0,\quad 
\partial_{p}f(x_*,p_*)\neq 0.
\ee 
Then a bifurcation occurs at $(x_*,p_*)$, which can be proven to be locally 
topologically equivalent to the fold (or saddle-node) bifurcation
normal form
\be
\label{eq:foldnf}
x'=p-x^2.
\ee
For $p>0$, there 
are two equilibrium points $x_\pm=\pm\sqrt{p}$. $x_+$ is locally
stable, while $x_-$ is unstable. For $p<0$, there are no equilibria. Only for 
$p>0$, there is a trapping region given by $\cV_{p>0}=[x_-,x_+]$.
For $p_{-} < 0$, $p_+ > 0$, we switch between 
$$ 
f_{-1}(x) = p_{-} - x^2, \quad f_1(x) = p_+ - x^2
$$ 
at rates $\lambda_{-}$ and $\lambda_+$.  

\begin{thm}     
\label{thm:divergence_fold}
Let $\nu$ be a probability measure on $\R \times I$, and let 
$$ 
\tau_a = \inf \{t \geq 0: X_t \leq a\}
$$ 
for $a < -\sqrt{p_+}$. Let $\nu \Pp^{\tau_a}$ be the law of $(X,E)$ with initial distribution $\nu$ and stopped at time $\tau_a$. Then, there is a nonincreasing function $g: (-\infty, -\sqrt{p_+}-2] \to (0,\infty)$ such that $\int_{-\infty}^{-\sqrt{p_+}-2} g(a) \ da < \infty$ and 
$$ 
\nu \Pp^{\tau_a}(\tau_{-\sqrt{p_+}-1} < \infty) = 1, \quad  \nu \Pp^{\tau_a}(\tau_a - \tau_{a+1} < g(a) \mid \tau_{-\sqrt{p_+}-1} < \infty) = 1, \quad a \leq -\sqrt{p_+} - 2. 
$$ 
\end{thm}

In words, $X_t$ diverges to $-\infty$ in finite time almost surely. 

\bpf[Proof of Theorem~\ref{thm:divergence_fold}] Let us show that $\nu \Pp^{\tau_a}(\tau_{-\sqrt{p_+}-1} < \infty) = 1$. The rest is analogous to the proof of Proposition~\ref{prop:divergence_to_infinity}. For fixed $x \geq -\sqrt{p_+}-1$, let $s_1 \geq 0$ such that $\Phi_{-1}^{s_1}(x) = -\sqrt{p_+}-1$, and let $s_2 > 0$ such that $\Phi_{-1}^{s_2}(\sqrt{p_+}) = -\sqrt{p_+}-1$.  Set $s = \max\{s_1, s_2\}$ and let $i \in I$.  With $\delta_{x,i} \Pp^{\tau_a}$-probability $1$, we have $\tau_{-\sqrt{p_+}-1} < \infty$ or there is $r \geq 0$ such that $E_t = -1$ for all $t \in [r,r+s]$. But the latter case also implies $\tau_{-\sqrt{p_+}-1} < \infty$ because any switching trajectory starting from $x$ cannot move to the right of $\max\{x, \sqrt{p_+}\}$. As a result, 
$$ 
\nu \Pp^{\tau_a} \left(\tau_{-\sqrt{p_+}-1} < \infty \right) = \sum_{i \in I} \int_{-\infty}^{\infty} \delta_{x,i} \Pp^{\tau_a}\left(\tau_{-\sqrt{p_+}-1} < \infty \right) \ \nu(dx \times \{i\}) = 1. 
$$     
\epf

\medskip

%%%%%%%%%%%%%%%%%%%%%%%%%%%%%%%%%%%%%%%%%%%%%%%%%%%%%%%%%%%%%%%%%%%%%%%%%%%%%%%%
\section{Applications}
\label{sec:applications}

In this section, we provide several very brief examples of systems where the 
switching viewpoint near bifurcations via PDMPs can yield insight into concrete 
dynamical systems arising in applications. In particular, the normal form 
results can be used sufficiently close to bifurcation points after a normal form
transformation. Furthermore, they can also be used directly to form conjectures 
about the dynamics of the applications.

%%%%%%%%%%%%%%%%%%%%%%%%%%%%%%%%%%%%%%%%%%%%%%%%%%%%%%%%%%%%%%%%%%%%%%%%%%%%%%%%
\subsection{The Paradox of Enrichment}
\label{ssec:enrichment}

The paradox of enrichment is a classical topic in ecology. One simple variant
can be found in classical predator-prey systems, such as the 
Rosenzweig-MacArthur~\cite{RosenzweigMacArthur} model
\be
\label{eq:RM}
\begin{array}{lcl}
x'&=& x\left(1-\frac{x}{p}\right)-\frac{xy}{1+x},\\
y'&=& \beta\frac{xy}{1+x} - my,
\end{array}
\ee
where $x,y\in[0,\I)$ are population densities of prey and predator, $\beta>0$ is a 
parameter representing a conversion factor, $m>0$ is the mortality of the 
predator, while $p>0$ is the carrying capacity for the prey. The basic 
concept of the paradox of enrichment~\cite{Rosenzweig} is that increasing 
the carrying capacity $p>0$ can actually lead to more likely extinction events 
triggered by additional stochastic effects, which can be supported by a classical
bifurcation analysis of~\eqref{eq:RM} as follows: Let us fix 
\benn
m=1 \qquad \text{and}\qquad \beta=3
\eenn
while varying $p$ as the main bifurcation parameter. Besides the two 
boundary equilibrium points $(x,y)=(0,0)$ and $(x,y)=(p,0)$, we find the
nontrivial co-existence equilibrium point 
\benn
(x_*,y_*)=(x_*,y_*(p))=\left(\frac12,\frac{3(2p-1)}{4p}\right)
\eenn
which is in the relevant domain given by the non-negative quadrant for $p>1/2$.
Linearizing~\eqref{eq:RM} around $(x_*,y_*)$ shows that the coexistence 
equilibrium is locally asymptotically stable for $p\in(0,2)$. Another direct
calculation shows that a supercritical Hopf bifurcation occurs at $p_*=2$.
The resulting locally asymptotically stable limit cycle generated in the Hopf 
bifurcation for $p>2$ grows in phase space. Hence, solutions can get closer to 
the two coordinate axes $\{x=0,y\geq 0\}$ and $\{y=0,x\geq 0\}$, which could 
make it more likely that a stochastic effect triggers an extinction event of
a species. Therefore, enrichment may lead to a potential increase in
extinction events. Of course, it is important to mention that there is still 
a debate in the literature on the mechanisms and possible variations of the 
paradox of enrichment~\cite{Kirk,KuehnSDEcont1,McCauleyMurdoch}. We do not 
provide here a full discussion of the various arguments made in favor or 
against the paradox but instead point towards the effect of randomness in
the parameter $p$.\medskip

From an ecological perspective, it can be plausible to view $p$ as a 
parameter, which switches randomly between different carrying capacities 
since the environment might be driven by 
external random events such as droughts, floods, storms, earthquakes, 
sudden human intervention, or even just different seasonal climate
conditions. Suppose we switch $p$ randomly in a range near $p_*=2$ with
rates $\lambda_\pm$ and values $p_\pm$ as defined in Section~\ref{ssec:superHopf},
where the parameters $p_\pm$ are chosen so that the supercritical Hopf 
normal form~\eqref{eq:superHopfnf} is a good local approximation 
of~\eqref{eq:RM} near $p_*$. Then Theorem~\ref{thm:existence_hopf}
suggests an interesting dichotomy of the ergodic IPM. Either, we have
only the invariant measure $\delta$, which is concentrated on the equilibrium 
branch $(x_*,y_*(p))$ with probabilities determined by the switching rates,
or we have two probability measures given by $\delta$ and $\nu$, where $\nu$
is a non-trivial product measure also supported on the periodic orbit. One
natural ecological interpretation of this effect is that we can actually 
avoid the paradox of enrichment from the viewpoint of measures if we restrict
to those switching rates, which only lead to the IPM $\delta$, i.e., that we
switch frequently enough from the periodic stable regime above the bifurcation 
to the stationary stable regime below the bifurcation.

%%%%%%%%%%%%%%%%%%%%%%%%%%%%%%%%%%%%%%%%%%%%%%%%%%%%%%%%%%%%%%%%%%%%%%%%%%%%%%%%
\subsection{Relaxation Oscillations}
\label{ssec:relaxation}

Consider the van der Pol (vdP) / FitzHugh-Nagumo (FHN) system
\benn
\begin{array}{lcl}
x' &=& p-\frac13 x^3+x,\\
p' &=& -\varepsilon x,
\end{array}
\eenn
which is a classical model used for bistable systems with relaxation
oscillations and excitability. The parameter $\varepsilon$ is usually assumed to be small, so
that in the singular limit $\varepsilon\ra 0$, we obtain the fast subsystem 
ODE
\be
\label{eq:vdPfss}
x'=p-\frac13 x^3+x,
\ee
where $p\in\R$ becomes a parameter. We can also view $p$ as a random parameter
for the dynamics. Observe that there are several branches of equilibrium points 
for~\eqref{eq:vdPfss} given by solving
\benn
p=\frac13 x^3-x.
\eenn
If $p\in(-2/3,2/3)$, then there are three equilibria, two locally asymptotically
stable and one unstable. At $p_*=-2/3$ and $p_*=2/3$, there are non-degenerate
fold bifurcations. While for $|p|>2/3$, there is always only one globally 
stable equilibrium point $x_*$. Let us focus on the case of switching $p$
near $p_*=2/3$; the case $p_*=-2/3$ simply follows by a symmetry argument.
If we switch the dynamics randomly above and below the fold bifurcation, 
Theorem~\ref{thm:divergence_fold} suggests that with probability one, we are
going to diverge away from the region of the fold, i.e., we are going to
obtain a point measure eventually concentrated on single remaining 
equilibrium $x_*$ existing for $p_*>2/3$. Hence, if we have random
switching across \emph{both} folds, then it is possible to obtain the
classical structure of a relaxation oscillations~\cite{MisRoz}.
This confirms similar observations made already numerically for a similar 
class of randomly switched van der Pol oscillators in~\cite{KuehnQuenched}.

%%%%%%%%%%%%%%%%%%%%%%%%%%%%%%%%%%%%%%%%%%%%%%%%%%%%%%%%%%%%%%%%%%%%%%%%%%%%%%%%
\subsection{Adaptive Swarming}
\label{ssec:swarming}

The next ODE model we are going to discuss is motivated by the swarming motion of
locusts in a ring-shaped arena~\cite{Buhletal}. An adaptive network~\cite{GrossSayama} 
model for this situation was proposed in~\cite{HuepeZschalerDoGross}. The network
model views locusts in clockwise-moving and anti-clockwise-moving nodes and keeps 
track of interactions between different locusts/nodes via the links of the network. The
model is reduced to a low-dimensional ODE via moment closure~\cite{KuehnMC} and
we focus here only on the essential features of the following low-dimensional
ODE model
\be
\label{eq:Huepe1}
\begin{array}{lcl}
x_1'&=& q(x_1-x_2)+w_3(y_3^2/(2x_2)-y_3^2/(2x_1)),\\
x_2'&=& q(x_2-x_1)+w_3(y_3^2/(2x_1)-y_3^2/(2x_2)),\\
y_1'&=& q(y_3-2y_1)+w_2(y_3+y_3^2/L-2y_3y_1/x_1),\\
&&+w_3(y_3^2/x_2+y_3^3/(2x_2^2)-y_3^2y_1/x_1^2)+ a_e x_1^2 - d_e y_1\\
y_2' &=& q(y_3-2y_2)+w_2(y_3+y_3^2/x_1-2y_3y_1/x_1),\\
&&+w_3(y_3^2/x_1+y_3^3/(2x_1^2)-y_3^2y_2/x_2^2)+a_e  x_2^2 - d_e y_2
\end{array}
\ee 
and the conservation equation
\be
\label{eq:Huepe2}
(y_1+y_2+y_3)'=a_0x_1x_2-d_0 y_3+a_e (x_1^2+x_2^2)-d_e(y_1+y_2),
\ee
where $q,w_2,w_3,a_e,d_e,a_0,d_0$ are positive parameters. Basically, $x_1$ and $x_2$ 
correspond to proportions of clockwise $R$ (``right'') and anti-clockwise 
$L$ (``left'') moving nodes, while $y_{1,2,3}$ capture the link densities $RR$, $LL$ 
and $RL$ between the two classes of nodes respectively. One checks that for $a_e=0=d_e$,
one can solve the steady-state problem, which provides a branch of solutions 
given by
\benn
x_1=\frac12=x_2.
\eenn
This state corresponds to an equal number of left and right moving nodes.
This disordered state is locally asymptotically stable up to a 
supercritical pitchfork bifurcation at $a_0^*=2d_0\sqrt{2q/w_3}$. The 
parameter $a_0$ controls the rate at which new connections between left
and right moving nodes form. Therefore, for a high connectivity between
different groups, the system can move into an ordered phase given by 
the steady state
\benn
(x_1)_\pm:=\frac12\pm\frac12\sqrt{1-8qd_0^2/(w_3a_0^2)}
\eenn
and similarly for $x_2$ with reversed signs. This corresponds to a classical
symmetry-breaking and above the supercritical pitchfork, the two majority
states are locally asymptotically stable.
Clearly, we can also view $a_0=:p$ as our randomly switched parameter across the
pitchfork bifurcation. Then Theorem~\ref{thm:existence_pitchfork} provides us
with the case of either one or three IPM if we switch near $a_0^*$. The 
interpretation for swarming is that we effectively can allow for a certain
percentage of disordered motion as long as the switching rate back into the ordered
phase is large enough to get an effective ordered phase. Furthermore, if we 
have the case of three IPM, then we are bound to observe not only a pure ordered
state but intermittent phases of disordered motion as the non-trivial measures are 
supported also near the locally unstable state above the bifurcation point. 

%%%%%%%%%%%%%%%%%%%%%%%%%%%%%%%%%%%%%%%%%%%%%%%%%%%%%%%%%%%%%%%%%%%%%%%%%%%%%%%%
%\section{Conclusion}
%\label{sec:conclusion}

%%%%%%%%%%%%%%%%%%%%%%%%%%%%%%%%%%%%%%%%%%%%%%%%%%%%%%%%%%%%%%%%%%%%%%%%%%%%%%%%
%%%%%%%%%%%%%%%%%%%%%%%%%%%%%%%%%%%%%%%%%%%%%%%%%%%%%%%%%%%%%%%%%%%%%%%%%%%%%%%%
%%%%%%%%%%%%%%%%%%%%%%%%%%%%%%%%%%%%%%%%%%%%%%%%%%%%%%%%%%%%%%%%%%%%%%%%%%%%%%%%

%\newpage 

\bibliographystyle{plain}
\bibliography{bif_switch}

\begin{thebibliography}{10}

\bibitem{Anisimov}
Vladimir~V. Anisimov.
\newblock {\em Switching Processes in Queueing Models}.
\newblock ISTE Ltd, John Wiley \& Sons, Inc., 2008.

\bibitem{ArnoldSDE}
L.~Arnold.
\newblock {\em Random Dynamical Systems}.
\newblock Springer, Berlin Heidelberg, Germany, 2003.

\bibitem{Bakhtin}
Yuri Bakhtin and Tobias Hurth.
\newblock Invariant densities for dynamical systems with random switching.
\newblock {\em Nonlinearity}, 25(10):2937--2952, 2012.

\bibitem{Mattingly}
Yuri Bakhtin, Tobias Hurth, and Jonathan~C. Mattingly.
\newblock Regularity of invariant densities for 1d-systems with random
  switching.
\newblock {\em Nonlinearity}, 28:3755--3787, 2015.

\bibitem{Balazs}
M.~Bal\'{a}zs, G.~Horv\'{a}th, S.~Kolumb\'{a}n, P.~Kov\'{a}cs, and M.~Telek.
\newblock Fluid level dependent {M}arkov fluid models with continuous zero
  transition.
\newblock {\em Performance Evaluation}, 68(11):1149 -- 1161, 2011.
\newblock Special Issue: Performance 2011.

\bibitem{BaudelBerglund}
M.~Baudel and N.~Berglund.
\newblock Spectral theory for random {Poincar\'e} maps.
\newblock {\em SIAM J. Math. Anal.}, 49(6):4319--4375, 2017.

\bibitem{B17}
M.~Bena{\"\i}m.
\newblock Stochastic persistence, part {I}.
\newblock 2017.
\newblock preprint.

\bibitem{BS17}
M.~Bena{\"\i}m and E.~Strickler.
\newblock Random switching between vector fields having a common zero.
\newblock 2017.
\newblock https://arxiv.org/abs/1702.03089.

\bibitem{Benaim}
Michel Bena{\"i}m, St{\'e}phane Le~Borgne, Florent Malrieu, and
  Pierre-Andr{\'e} Zitt.
\newblock Qualitative properties of certain piecewise deterministic {M}arkov
  processes.
\newblock {\em Annales de l'Institut Henri Poincar\'{e}}, 51:1040 -- 1075,
  2012.

\bibitem{Le_Borgne}
Michel Bena{\"{\i}}m, St{\'e}phane Le~Borgne, Florent Malrieu, and
  Pierre-Andr{\'e} Zitt.
\newblock Quantitative ergodicity for some switched dynamical systems.
\newblock {\em Electron. Commun. Probab.}, 17:no. 56, 14, 2012.

\bibitem{Zitt}
Michel Bena{\"{\i}}m, St{\'e}phane Le~Borgne, Florent Malrieu, and
  Pierre-Andr{\'e} Zitt.
\newblock On the stability of planar randomly switched systems.
\newblock {\em Ann. Appl. Probab.}, 24(1):292--311, 2014.

\bibitem{Lobry}
Michel Benaïm and Claude Lobry.
\newblock Lotka–volterra with randomly fluctuating environments or “how
  switching between beneficial environments can make survival harder”.
\newblock {\em Ann. Appl. Probab.}, 26(6):3754--3785, 12 2016.

\bibitem{BerglundGentz}
N.~Berglund and B.~Gentz.
\newblock {\em Noise-Induced Phenomena in Slow-Fast Dynamical Systems}.
\newblock Springer, 2006.

\bibitem{BredenKuehn}
M.~Breden and C.~Kuehn.
\newblock Exploring invariant sets of random dynamical systems via polynomial
  chaos.
\newblock {\em preprint}, pages 1--, 2019.

\bibitem{Bressloff}
Paul~C Bressloff.
\newblock Stochastic switching in biology: from genotype to phenotype.
\newblock {\em Journal of Physics A: Mathematical and Theoretical},
  50(13):133001, 2017.

\bibitem{Buhletal}
J.~Buhl, D.J. Sumpter, I.D. Couzin, J.J. Hale, E.~Despland, E.R. Miller, and
  S.J. Simpson.
\newblock From disorder to order in marching locusts.
\newblock {\em Science}, 312(5778):1402--1406, 2006.

\bibitem{DaPrato}
G.~Da~Prato and J.~Zabczyk.
\newblock {\em Ergodicity for infinite-dimensional systems}, volume 229 of {\em
  London Mathematical Society Lecture Note Series}.
\newblock Cambridge University Press, Cambridge, 1996.

\bibitem{Davis_article}
M.~H.~A. Davis.
\newblock Piecewise-deterministic {M}arkov processes: a general class of
  nondiffusion stochastic models.
\newblock {\em J. Roy. Statist. Soc. Ser. B}, 46(3):353--388, 1984.
\newblock With discussion.

\bibitem{Davis}
M.~H.~A. Davis.
\newblock {\em Markov models and optimization}, volume~49 of {\em Monographs on
  Statistics and Applied Probability}.
\newblock Chapman \& Hall, London, 1993.

\bibitem{EngelLambRasmussen}
M.~Engel, J.S. Lamb, and M.~Rasmussen.
\newblock Bifurcation analysis of a stochastically driven limit cycle.
\newblock {\em arXiv:1606.01137}, pages 1--, 2016.

\bibitem{Ethier_Kurtz}
Stewart~N. Ethier and Thomas~G. Kurtz.
\newblock {\em {M}arkov Processes}.
\newblock Wiley Series in Probability and Statistics. John Wiley \& Sons, Inc.,
  Hoboken, New Jersey, 2005.
\newblock Characterization and Convergence.

\bibitem{Faggionato}
A.~Faggionato, D.~Gabrielli, and M.~Ribezzi~Crivellari.
\newblock Non-equilibrium thermodynamics of piecewise deterministic {M}arkov
  processes.
\newblock {\em J. Stat. Phys.}, 137(2):259--304, 2009.

\bibitem{Gabrielli}
A.~Faggionato, D.~Gabrielli, and M.~Ribezzi~Crivellari.
\newblock Averaging and large deviation principles for fully-coupled piecewise
  deterministic {M}arkov processes and applications to molecular motors.
\newblock {\em {M}arkov processes and related fields}, 16:497--548, 2010.

\bibitem{Goldstein}
S.~Goldstein.
\newblock On diffusion by discontinuous movements, and on the telegraph
  equation.
\newblock {\em Quart. J. Mech. Appl. Math.}, 4:129--156, 1951.

\bibitem{GrossSayama}
T.~Gross and H.~Sayama, editors.
\newblock {\em Adaptive Networks: Theory, Models and Applications}.
\newblock Springer, 2009.

\bibitem{GH}
J.~Guckenheimer and P.~Holmes.
\newblock {\em Nonlinear Oscillations, Dynamical Systems, and Bifurcations of
  Vector Fields}.
\newblock Springer, New York, NY, 1983.

\bibitem{Hairer:ergodicity-lectures}
Martin Hairer.
\newblock Ergodic properties of {M}arkov processes.
\newblock Lectures given at the University of Warwick,
  http://www.hairer.org/notes/Markov.pdf, 2006.

\bibitem{Hersh}
Reuben Hersh.
\newblock The birth of random evolutions.
\newblock {\em Math. Intelligencer}, 25(1):53--60, 2003.

\bibitem{HuepeZschalerDoGross}
C.~Huepe, G.~Zschaler, A.-L. Do, and T.~Gross.
\newblock Adaptive network models of swarm dynamics.
\newblock {\em New J. Phys.}, page (073022), 2011.

\bibitem{Kac}
Mark Kac.
\newblock A stochastic model related to the telegrapher's equation.
\newblock {\em Rocky Mountain J. Math.}, 4:497--509, 1974.

\bibitem{Kirk}
K.L. Kirk.
\newblock Enrichment can stabilize population dynamics: autotoxins and density
  dependence.
\newblock {\em Ecol.}, 79:2456--2462, 1998.

\bibitem{KuehnSDEcont1}
C.~Kuehn.
\newblock Deterministic continuation of stochastic metastable equilibria via
  {Lyapunov} equations and ellipsoids.
\newblock {\em SIAM J. Sci. Comp.}, 34(3):A1635--A1658, 2012.

\bibitem{KuehnMC}
C.~Kuehn.
\newblock Moment closure - a brief review.
\newblock In E.~Sch{\"o}ll, S.~Klapp, and P.~H{\"o}vel, editors, {\em Control
  of Self-Organizing Nonlinear Systems}, pages 253--271. Springer, 2016.

\bibitem{KuehnQuenched}
C.~Kuehn.
\newblock Quenched noise and nonlinear oscillations in bistable multiscale
  systems.
\newblock {\em EPL (Europhysics Letters)}, 120:10001, 2017.

\bibitem{Kuznetsov}
Yu.A. Kuznetsov.
\newblock {\em Elements of Applied Bifurcation Theory}.
\newblock Springer, New York, NY, 3rd edition, 2004.

\bibitem{Lawley}
Sean~D. Lawley, Jonathan~C. Mattingly, and Michael~C. Reed.
\newblock Sensitivity to switching rates in stochastically switched {ODE}s.
\newblock {\em Commun. Math. Sci.}, 12(7):1343--1352, 2014.

\bibitem{Cui}
Dan Li, Shengqiang Liu, and Jing'an Cui.
\newblock Threshold dynamics and ergodicity of an sirs epidemic model with
  markovian switching.
\newblock {\em Journal of Differential Equations}, 263(12):8873 -- 8915, 2017.

\bibitem{Majda}
Andrew~J. Majda and Xin~T. Tong.
\newblock Geometric ergodicity for piecewise contracting processes with
  applications for tropical stochastic lattice models.
\newblock {\em Communications on Pure and Applied Mathematics},
  69(6):1110--1153, 2016.

\bibitem{Malrieu_2015}
Florent Malrieu.
\newblock Some simple but challenging {M}arkov processes.
\newblock {\em Ann. Fac. Sci. Toulouse Math. (6)}, 24(4):857--883, 2015.

\bibitem{McCauleyMurdoch}
E.~McCauley and W.W. Murdoch.
\newblock Predator-prey dynamics in environments rich and poor in nutrients.
\newblock {\em Nature}, 343:455--461, 1990.

\bibitem{MisRoz}
E.F. Mishchenko and N.Kh. Rozov.
\newblock {\em Differential Equations with Small Parameters and Relaxation
  Oscillations (translated from Russian)}.
\newblock Plenum Press, 1980.

\bibitem{NairSarkarSujith}
V.~Nair, S.~Sarkar, and R.I. Sujith.
\newblock Uncertainty quantification of subcritical bifurcations.
\newblock {\em Probab. Eng. Mech.}, 34:177--188, 2013.

\bibitem{Nualart}
David Nualart.
\newblock {\em The {M}alliavin calculus and related topics}.
\newblock Probability and its Applications (New York). Springer-Verlag, Berlin,
  second edition, 2006.

\bibitem{Rosenzweig}
M.L. Rosenzweig.
\newblock Paradox of enrichment: Destabilization of exploitation ecosystems in
  ecological time.
\newblock {\em Science}, 171:385--387, 1971.

\bibitem{RosenzweigMacArthur}
M.L. Rosenzweig and R.H. MacArthur.
\newblock Graphical representation and stability conditions of predator-prey
  interactions.
\newblock {\em American Naturalist}, 97:209--223, 1963.

\bibitem{SadhuKuehn}
S.~Sadhu and C.~Kuehn.
\newblock Stochastic mixed-mode oscillations in a three-species predator-prey
  model.
\newblock {\em Chaos}, 28(3):033606, 2017.

\bibitem{Strickler}
Edouard Strickler.
\newblock Randomly switched vector fields sharing a zero on a common invariant
  face.
\newblock {\em Available at https://arxiv.org/abs/1810.06331}, 2018.

\bibitem{Strogatz}
S.H. Strogatz.
\newblock {\em Nonlinear Dynamics and Chaos}.
\newblock Westview Press, 2000.

\bibitem{Teschl}
G.~Teschl.
\newblock {\em Ordinary Differential Equations and Dynamical Systems}.
\newblock AMS, 2012.

\bibitem{Yin}
G.~George Yin and Chao Zhu.
\newblock {\em Hybrid switching diffusions}, volume~63 of {\em Stochastic
  Modelling and Applied Probability}.
\newblock Springer, New York, 2010.
\newblock Properties and applications.

\end{thebibliography}
%\bibliography{../my_refs}

\end{document}